\documentclass[11pt, twoside]{article}

\usepackage[english]{babel}

\usepackage{bbm}
\usepackage{amssymb}
\usepackage{amsfonts}
\usepackage{amsmath}
\usepackage{amsthm}
\usepackage{color}
\usepackage{mathrsfs}
\usepackage{txfonts}
\usepackage{enumerate}
\usepackage{anysize}
\usepackage{indentfirst}

\usepackage{latexsym}

\usepackage[colorlinks=true,
linkcolor=blue,
citecolor=red,
urlcolor=magenta,
]{hyperref}

\allowdisplaybreaks

\newtheorem{theorem}{Theorem}[section]
\newtheorem{lemma}[theorem]{Lemma}
\newtheorem{corollary}[theorem]{Corollary}
\newtheorem{proposition}[theorem]{Proposition}

\theoremstyle{definition}
\newtheorem{remark}[theorem]{Remark}

\numberwithin{equation}{section}

\begin{document}

\arraycolsep=1pt

\title{\bf\Large Parabolic Extrapolation of Compactness with Applications
to Characterizations of the Parabolic CMO Space via 
Parabolic Fractional Commutators
\footnotetext{\hspace{-0.35cm} 2020
\emph{Mathematics Subject Classification}.
Primary 47B47; Secondary 42B25, 42B35, 46B70, 47A30.
\endgraf \emph{Key words and phrases}.
extrapolation and interpolation of compactness,
parabolic fractional integral,
parabolic fractional commutator,
parabolic Muckenhoupt weight with time lag,
parabolic CMO space.
\endgraf
This project is partially supported by the National Natural
Science Foundation of China (Grant Nos. 12431006, 12371093, and 12501129), the Beijing
Natural Science Foundation (Grant No. 1262011),
and the Fundamental Research Funds for the Central Universities
(Grant No. 2233300008).}}
\author{Mingming Cao, Weiyi Kong, Dachun Yang,
Wen Yuan and Chenfeng Zhu}
\maketitle

\vspace{-0.8cm}

\begin{center}
\begin{minipage}{13cm}
{\small {\bf Abstract}\quad In this article, we prove two compactness parabolic
extrapolation theorems for linear operators and their commutators on
parabolic weighted
Lebesgue spaces. As applications, we show that $b\in\mathrm{PCMO}$ (parabolic
CMO space) if and only if $b\in L_\mathrm{loc}^1$ and the commutator
$[b,I_\alpha^{\gamma+}]_k$ is compact from $L^r(\omega^r)$ to $L^q(\omega^q)$
for any $k\in\mathbb{N}$, $1<r<q<\infty$ with
$\frac{1}{r}-\frac{1}{q}=\alpha\in(0,1)$, and $\omega\in A_{r,q}^+(\gamma)$
(parabolic Muckenhoupt class), where $I_\alpha^{\gamma+}$ is the parabolic
fractional integral. Several other equivalent characterizations of
the space $\mathrm{PCMO}$ are also given. Our approach combines the unweighted
compactness for mixed commutators, a local lower estimate
for $[b,I_\alpha^{\gamma+}]_k$, and both a
quantitative Uchiyama-type description and the equivalence with
$\inf_T\{\|[\cdot,I_\alpha^{\gamma+}]_k-T\|_{L^r\to L^q}^\frac{1}{k}\}$ of
$\mathrm{dist}_{\mathrm{PBMO}}(\cdot,\mathrm{PCMO})$, where the infimum is
taken over all compact operators $T:L^r\to L^q$. Finally, we establish two
parabolic weighted Fr\'echet--Kolmogorov compactness criteria, which further
induce a parabolic interpolation theorem of compactness.
}
\end{minipage}
\end{center}

\vspace{0.2cm}

\tableofcontents

\vspace{0.2cm}

\section{Introduction}

Throughout this article, we work in $\mathbb{R}^{n+1}$ and,
unless otherwise specified, we take $\mathbb{R}^{n+1}$ as the default
underlying space. Muckenhoupt weights form a basic framework for weighted norm
inequalities in harmonic analysis and partial differential equations.
Parabolic Muckenhoupt weights with time lag, introduced and developed by
Kinnunen and Saari \cite{ks(na-2016), ks(apde-2016)} and Kinnunen and
Myyryl\"ainen \cite{km(am-2024),2310.00370}, are motivated by two complementary
considerations. On the one hand, these Parabolic Muckenhoupt weights
with time lag provide an appropriate
$(n+1)$-dimensional form of one-dimensional one-sided weights introduced by
Sawyer \cite{s(tams-1986)}. On the other hand, they are closely related to the
regularity theory of the doubly nonlinear parabolic equation
\begin{align}\label{20261012.2050}
\partial_t\left(|u|^{p-2}u\right)
-\mathrm{div}\left(|\nabla u|^{p-2}\nabla u\right)=0.
\end{align}
Here, and thereafter, we \emph{always fix} $p\in(1,\infty)$. Let
$(x,t)\in\mathbb{R}^{n+1}$ and $L\in(0,\infty)$. A \emph{parabolic rectangle}
$R(x,t,L)$ centered at $(x,t)$ with edge length $l(R):=L$ is defined by setting
\begin{align*}
R:=R(x,t,L):=Q(x,L)\times[t-L^p,t+L^p),
\end{align*}
where $Q(x,L):=\{y:=(y_1,\dots,y_n)\in\mathbb{R}^n:x_i-L\leq
y_i<x_i+L\mbox{\ for\ any\ }i\in\mathbb{N}\cap[1,n]\}$
is the \emph{cube} in $\mathbb{R}^n$ centered at $x$ with edge length $2L$.
Moreover, for any given $\gamma\in[0,1)$, the \emph{$\gamma$-lower part
$R^-(\gamma)$} and the \emph{$\gamma$-upper part $R^+(\gamma)$} of $R$ are
defined, respectively, by setting
\begin{align*}
R^-(\gamma):=Q(x,L)\times\left[t-L^p,t-\gamma L^p\right)
\end{align*}
and
\begin{align*}
R^+(\gamma):=Q(x,L)\times\left[t+\gamma L^p,t+L^p\right),
\end{align*}
where $\gamma$ is called the \emph{time lag}. Let $\mathcal{R}$ be the set of
all parabolic rectangles in $\mathbb{R}^{n+1}$. The parabolic scaling
$(x,t)\mapsto(\lambda x,\lambda^pt)$ with $\lambda\in(0,\infty)$ is inspired by
estimates concerning the regularity theory of \eqref{20261012.2050}, and the time lag
$\gamma$ appears naturally in the Harnack inequality associated with
\eqref{20261012.2050} (see, for example, \cite{kk(ma-2007), s(rmi-2016)}).

By a \emph{weight}, we mean a nonnegative locally integrable function on
$\mathbb{R}^{n+1}$ which is positive almost everywhere. In what follows, if there is
no confusion, then we \emph{always omit} the differential $dx\,dt$ in all integral
representations to simplify the presentation. For any measurable subset
$A\subset\mathbb{R}^{n+1}$ with $|A|\in(0,\infty)$ and for
any integrable function $f$ on $A$, let
\begin{align*}
f_A:=\fint_Af:=\frac{1}{|A|}\int_Af.
\end{align*}
Let $\gamma\in[0,1)$ and $1\leq r\leq q<\infty$. The \emph{parabolic
Muckenhoupt class $A_{r,q}^+(\gamma)$ with time lag $\gamma$} is defined to be
the set of all weights $\omega$ on $\mathbb{R}^{n+1}$ such that
\begin{align*}
[\omega]_{A_{r,q}^+(\gamma)}:=
\sup_{R\in\mathcal{R}}\left[\fint_{R^-(\gamma)}\omega^q\right]^\frac1q
\left[\fint_{R^+(\gamma)}\omega^{-r'}\right]^\frac{1}{r'}<\infty
\end{align*}
with the usual modification made when $r=1$ (see, for instance,
\cite[(1.2)]{mhy(fm-2023)}). If the condition above holds with
the time axis reversed, then $\omega\in A_{r,q}^-(\gamma)$ (the other
\emph{parabolic Muckenhoupt class with time lag $\gamma$}).
We denote by $A_q^+(\gamma)$ the set of all weights $\omega$ on $\mathbb{R}^{n+1}$
such that $\omega^\frac{1}{q}\in A_{q,q}^+(\gamma)$, and
let $[\omega]_{A_q^+(\gamma)}:=[\omega^\frac{1}{q}]_{A_{q,q}^+(\gamma)}^q$
(see, for example, \cite{km(am-2024),ks(apde-2016)}). Let
$A_\infty^+(\gamma):=\bigcup_{q\in[1,\infty)}A_{q}^+(\gamma)$.
The classes $A_q^-(\gamma)$ and $A_\infty^-(\gamma)$
are defined similarly. It is worth noting that any positive solution to
\eqref{20261012.2050} is a parabolic Muckenhoupt weight with time lag
(see, for instance, \cite[Proposition 8.3]{ks(apde-2016)} and \cite[Proposition
5.17]{kyy(arXiv:2602.09741)}). We refer to \cite{cm24, 2310.00370,
kmyz(pa-2023), kyy(arXiv:2608.13307), kyyz(cjm-2025), lmv(arXiv:2604.12561),
mhy(fm-2023)} for more studies on parabolic Muckenhoupt weights
with time lag and to \cite{kmy(ma-2023), kyy(arXiv:2602.09741),
kyyz(cvpde-2025), mos(arXiv:2509.24486), my(mz-2024), s(rmi-2016),
s(ampa-2018)} for further research about related parabolic function spaces.

The extrapolation theorem, initially established by Rubio de Francia
\cite{r(ajm-1984)}, is one of the most useful tools in the study of weighted
norm inequalities. One-dimensional one-sided extrapolation was developed in
\cite{clm(jam-2018), lr(mia-2007), mr(prse-2000)}. In our previous work
\cite{ckyyz(jga-2026)}, we gave two parabolic extrapolation theorems
of boundedness and applied them to characterize $\mathrm{PBMO}$ (parabolic BMO space)
through the weighted boundedness of the first-order parabolic fractional commutators.
For classical weights, Hyt\"onen and Lappas \cite{hl(rmi-2023)}
established the extrapolation of compactness for linear operators.
Their compactness extrapolation theorems have recently been
extended to various settings and frameworks (see, for example,
\cite{cdlsy(arXiv:2026.01523), coy(tams-2022), hl(im-2022),
kl(arXiv:2512.09721), kl(arXiv:2605.27165), l(sm-2022), lwy(cm-2023),
ln(jfaa-2024)}). Mart\'in-Reyes and Rivera-R\'ios \cite{2406.11663v2}
proved one-dimensional one-sided counterparts of the compactness
extrapolation theorems.

Recall that the \emph{parabolic distance} $d_p$ on
$\mathbb{R}^{n+1}\times\mathbb{R}^{n+1}$ is defined by setting,
for any $(x,t):=(x_1,\dots,x_n,t),(y,s):=(y_1,\dots,y_n,s)\in\mathbb{R}^{n+1}$,
\begin{align*}
d_p((x,t),(y,s)):=\max\left\{\|x-y\|_\infty,\,|t-s|^{\frac{1}{p}}\right\},
\end{align*}
where $\|x-y\|_\infty:=\max\{|x_1-y_1|,\dots,|x_n-y_n|\}$.
Also recall that, for any $q\in(0,\infty)$ and any weight $\omega$,
the \emph{weighted Lebesgue space $L^q(\omega)$} is defined to be the set of all
measurable functions $f$ on $\mathbb{R}^{n+1}$ such that
\begin{align*}
\|f\|_{L^q(\omega)}:=
\left(\int_{\mathbb{R}^{n+1}}|f|^q\omega\right)^\frac1q<\infty.
\end{align*}
If $\omega\equiv1$, then $L^q(\omega)$ reduces to the \emph{Lebesgue space
$L^q$}. Moreover, the \emph{Lebesgue space $L^\infty$} is defined to be the set
of all measurable functions $f$ on $\mathbb{R}^{n+1}$ such that
\begin{align*}
\|f\|_{L^\infty}:=\mathop\mathrm{ess\,sup}_{(x,t)\in\mathbb{R}^{n+1}}
|f(x,t)|<\infty.
\end{align*}
The following parabolic extrapolation
of compactness is the first main result of this article.

\begin{theorem}\label{parabolic extrapolation of compactness}
Let $\gamma\in(0,1)$, $1<r_0\leq q_0<\infty$, and $T$ be a linear operator
bounded from $L^{r_0}(\omega_0^{r_0})$ to $L^{q_0}(\omega_0^{q_0})$ for all
$\omega_0\in A_{r_0,q_0}^+(\gamma)$,
with the operator norm bounded by a constant
depending only on $n,p,\gamma,r_0,q_0$, and
$[\omega_0]_{A^+_{r_0,q_0}(\gamma)}$. If there exist $1<r_1\leq
q_1<\infty$ satisfying $\frac{1}{r_1}-\frac{1}{q_1}=
\frac{1}{r_0}-\frac{1}{q_0}$ such that $T$ is compact from
$L^{r_1}(\omega_1^{r_1})$ to $L^{q_1}(\omega_1^{q_1})$
for some $\omega_1\in A_{r_1,q_1}^-(\gamma)$, then,
for any $1<r\leq q<\infty$ with $\frac{1}{r}-\frac{1}{q}=
\frac{1}{r_0}-\frac{1}{q_0}$ and for any $\omega\in A_{r,q}^+(\gamma)$,
$T$ is compact from $L^r(\omega^r)$ to $L^q(\omega^q)$.
\end{theorem}

\begin{remark}
Let all the notation be as in Theorem \ref{parabolic extrapolation of compactness}.
We point out that the weight assumptions in
Theorem \ref{parabolic extrapolation of compactness} are natural,
and the distinction between the forward and backward Muckenhoupt
classes is essential. Indeed, the extrapolation of compactness
from a single weight in $A^-_{r_1,q_1}(\gamma)$ to all weights in
$A^+_{r,q}(\gamma)$ is a parabolic counterpart of the one-sided
phenomenon first established in \cite[\text{Theorem 1.4}]{2406.11663v2};
see also \cite[Corollaries 1.5 and 1.6]{2406.11663v2}.
As shown in Remark \ref{rmk:20260729.1620}, the reverse weight class
in the compactness assumption cannot, in general, be replaced
by the forward weight class. Also, recall that the
\emph{Muckenhoupt class $A_{r,q}$ associated with $d_p$}
is defined to be the set of all weights $\omega$ on $\mathbb{R}^{n+1}$ such that
\begin{align*}
[\omega]_{A_{r,q}}:=\sup_{R\in\mathcal{R}}
\left(\fint_R\omega^q\right)^\frac1q
\left(\fint_R\omega^{-r'}\right)^\frac{1}{r'}<\infty.
\end{align*}
Since $R^\pm(\gamma)\subset R$, it follows that $A_{r_1,q_1}\subset
A^+_{r_1,q_1}(\gamma)\cap A^-_{r_1,q_1}(\gamma)$, and hence the assumptions
on weights in Theorem \ref{parabolic extrapolation of compactness}
may be specialized to $\omega_1\in A_{r_1,q_1}$. In particular, we can take
$\omega_1\equiv1$, in which case the compactness assumption is exactly that $T$
is compact from $L^{r_1}$ to $L^{q_1}$.
\end{remark}

The key step in proving Theorem
\ref{parabolic extrapolation of compactness} is to represent the desired
weighted space $L^r(\omega^r)$ [resp. $L^q(\omega^q)$] as the interpolation
space between $L^{r_1}(\omega_1^{r_1})$ [resp. $L^{q_1}(\omega_1^{q_1})$] and
another specially constructed parabolic weighted Lebesgue space.
As a direct corollary of Theorem
\ref{parabolic extrapolation of compactness},
we obtain an extrapolation result concerning the compactness of commutators.
To present this result, we first introduce some
necessary notation. For any $q\in[1,\infty)$, denote by the
\emph{notation $L_\mathrm{loc}^q$} the set of all
locally $q$-integrable functions on $\mathbb{R}^{n+1}$.  The parabolic
\emph{$\mathrm{BMO}$ space $\mathrm{PBMO}$ (associated with $d_p$)} is
defined to be the space of all $f\in L_\mathrm{loc}^1$ such that
\begin{align*}
\|f\|_{\mathrm{PBMO}}:=\sup_{R\in\mathcal{R}}\fint_R|f-f_R|<\infty.
\end{align*}
Denote by $\mathscr{M}$ the set of all measurable functions on
$\mathbb{R}^{n+1}$. Let $T$ be a linear operator defined on a linear subspace
$D(T)$ of $L_\mathrm{loc}^1$ and taking values in $\mathscr{M}$, and let $b\in
L_\mathrm{loc}^1$. Recall that the \emph{commutator $[b,T]$ of $T$ with the
symbol $b$} is defined by setting, for any suitable function $f$ on
$\mathbb{R}^{n+1}$ and for any $(x,t)\in\mathbb{R}^{n+1}$,
\begin{align*}
[b,T]f(x,t):=[b,T]_1f(x,t):=b(x,t)Tf(x,t)-T(bf)(x,t).
\end{align*}
For any $k\in\mathbb{N}\cap[2,\infty)$, the \emph{$k$-th order commutator
$[b,T]_k$ of $T$ with the symbol $b$} is defined by the recursive formula
$[b,T]_k:=[b,[b,T]_{k-1}]$. We can easily verify that, for any given
$k\in\mathbb{N}$, any suitable function $f$ on $\mathbb{R}^{n+1}$, and any
$(x,t)\in\mathbb{R}^{n+1}$,
\begin{align*}
[b,T]_kf(x,t)=T\left([b(x,t)-b(\cdot)]^kf(\cdot)\right)(x,t).
\end{align*}
The following corollary follows immediately from Theorem
\ref{parabolic extrapolation of compactness} and the parabolic extrapolation for commutators
(see \cite[Theorem~3.6]{ckyyz(jga-2026)}).

\begin{corollary}\label{ExtraCC}
Let $\gamma\in(0,1)$, $1<r_0\leq q_0<\infty$, $k\in\mathbb{N}$,
$b\in\mathrm{PBMO}$, and $T:L_{\rm loc}^1\to\mathscr{M}$ be a linear operator.
Assume that, for any $\omega_0\in A_{r_0,q_0}^+(\gamma)$,
$T$ is bounded from $L^{r_0}(\omega_0^{r_0})$ to $L^{q_0}(\omega_0^{q_0})$
with its operator norm bounded by a constant
depending only on $n,p,\gamma,r_0,q_0$, and
$[\omega_0]_{A^+_{r_0,q_0}(\gamma)}$. If there exist $1<r_1\leq
q_1<\infty$ satisfying $\frac{1}{r_1}-\frac{1}{q_1}=
\frac{1}{r_0}-\frac{1}{q_0}$ such that $[b,T]_k$ is compact from
$L^{r_1}(\omega_1^{r_1})$ to $L^{q_1}(\omega_1^{q_1})$
for some $\omega_1\in A_{r_1,q_1}^-(\gamma)$, then,
for any $1<r\leq q<\infty$ with $\frac{1}{r}-\frac{1}{q}=
\frac{1}{r_0}-\frac{1}{q_0}$ and for any $\omega\in A_{r,q}^+(\gamma)$,
$[b,T]_k$ is compact from $L^r(\omega^r)$ to $L^q(\omega^q)$.
\end{corollary}

As an application of Corollary~\ref{ExtraCC},
we characterize the parabolic CMO space in terms of the
weighted compactness of parabolic fractional commutators.
Denote by $C_\mathrm{c}^\infty$ the space of all infinitely differentiable
functions on $\mathbb{R}^{n+1}$ with compact support. The \emph{parabolic $\mathrm{CMO}$
space $\mathrm{PCMO}$} (associated with the parabolic distance $d_p$)
is defined to be the closure of
$C_\mathrm{c}^\infty$ in the $\mathrm{PBMO}$ topology, i.e.,
$\mathrm{PCMO}:=\overline{C_\mathrm{c}^\infty}^{\|\cdot\|_{\mathrm{PBMO}}}$.
For any $b\in\mathrm{PBMO}$, the distance from $b$ to $\mathrm{PCMO}$ is
defined by setting
\begin{align*}
\mathrm{dist}_\mathrm{PBMO}(b,\mathrm{PCMO})
:=\inf\{\|b-a\|_\mathrm{PBMO}:a\in\mathrm{PCMO}\}.
\end{align*}
Let $\gamma,\alpha\in(0,1)$. The \emph{parabolic fractional integral operator
$I_\alpha^{\gamma+}$ with time lag $\gamma$}, initially introduced in our
previous work \cite[Definition 5.1]{kyyz(cjm-2025)}, is defined by setting, for
any suitable function $f$ on $\mathbb{R}^{n+1}$ and for
any $(x,t)\in\mathbb{R}^{n+1}$,
\begin{align*}
I_\alpha^{\gamma+}f(x,t):=\int_{\bigcup\limits_{L\in(0,\infty)}R(x,t,L)^+(\gamma)}
\frac{f(y,s)}{[d_p((x,t),(y,s))]^{(n+p)(1-\alpha)}}\,dy\,ds.
\end{align*}
Note that the above integration domain in $I_\alpha^{\gamma+}f(x,t)$ is given by
\begin{align*}
\bigcup_{L\in(0,\infty)}R(x,t,L)^+(\gamma)=\left\{(y,s)\in\mathbb{R}^{n+1}:
s-t>\gamma\|x-y\|_\infty^p\right\}=:\Omega_{(x,t)}^{\gamma+}.
\end{align*}
With these definitions, the \emph{parabolic fractional commutator
$[b,I_\alpha^{\gamma+}]_k$ with time lag $\gamma$}, introduced in our previous
work \cite[Definition 1.8]{ckyyz(jga-2026)}, is defined by setting, for any
suitable function $f$ on $\mathbb{R}^{n+1}$ and for any $(x,t)\in\mathbb{R}^{n+1}$,
\begin{align*}
\left[b,I_\alpha^{\gamma+}\right]_kf(x,t)
:=\int_{\bigcup_{L\in(0,\infty)}R(x,t,L)^+(\gamma)}
\frac{[b(x,t)-b(y,s)]^kf(y,s)}{[d_p((x,t),(y,s))]^{(n+p)(1-\alpha)}}\,dy\,ds.
\end{align*}
We then present several characterizations of PCMO via the weighted compactness
of $[b,I_\alpha^{\gamma+}]$, which is the second main result of this article.

\begin{theorem}\label{thm:20260714.2255}
Let $\gamma,\alpha\in(0,1)$ and $b:\mathbb{R}^{n+1}\to\mathbb{R}$.
Then the following statements are mutually equivalent.
\begin{enumerate}
\item[\rm(i)] $b\in\mathrm{PCMO}$.

\item[\rm(ii)] $b\in L_\mathrm{loc}^1$ and, for any $1<r<q<\infty$ satisfying
$\frac1r-\frac1q=\alpha$, any $k\in\mathbb{N}$, and any $\omega\in A_{r,q}^+(\gamma)$,
$[b,I_\alpha^{\gamma+}]_k$ is compact from $L^r(\omega^r)$ to $L^q(\omega^q)$.

\item[\rm(iii)] $b\in L_\mathrm{loc}^1$ and there exist $1<r_0<q_0<\infty$, satisfying
$\frac{1}{r_0}-\frac{1}{q_0}=\alpha$, and $k_0\in\mathbb{N}$ such that,
for any $\omega\in A_{r_0,q_0}^+(\gamma)$, $[b,I_\alpha^{\gamma+}]_{k_0}$ is
compact from $L^{r_0}(\omega^{r_0})$ to $L^{q_0}(\omega^{q_0})$.

\item[\rm(iv)] $b\in L_\mathrm{loc}^1$ and there exist $1<r_0<q_0<\infty$, satisfying
$\frac{1}{r_0}-\frac{1}{q_0}=\alpha$, and $k_0\in\mathbb{N}$ such that
$[b,I_\alpha^{\gamma+}]_{k_0}$ is compact from $L^{r_0}$ to $L^{q_0}$.

\item[\rm(v)] $b\in L_\mathrm{loc}^1$ and there exist $1<r_0<q_0<\infty$ satisfying
$\frac{1}{r_0}-\frac{1}{q_0}=\alpha$, $k_0\in\mathbb{N}$, and $\omega_0\in
A_{r_0,q_0}^-(\gamma)$ such that $[b,I_\alpha^{\gamma+}]_{k_0}$ is compact from
$L^{r_0}(\omega_0^{r_0})$ to $L^{q_0}(\omega_0^{q_0})$.
\end{enumerate}
\end{theorem}

\begin{remark}
Let all the notation be as in Theorem \ref{thm:20260714.2255}. We point out
that the weight from $A_{r_0,q_0}^-(\gamma)$ in Theorem \ref{thm:20260714.2255}(v)
cannot be replaced by the weight from $A^+_{r_0,q_0}(\gamma)$.
Indeed, the boundedness with respect to a single parabolic Muckenhoupt weight need
not even imply that the symbol belongs to $\mathrm{PBMO}$. More precisely,
there exist $b_1\in L^1_{\mathrm{loc}}\setminus\mathrm{PBMO}$,
$b_2\in\mathrm{PBMO}\setminus\mathrm{PCMO}$, and a parabolic weight
$\omega_0$ such that, for any $1<r<q<\infty$ with $\frac1r-\frac1q=\alpha$ and
for any $k\in\mathbb{N}$, $\omega_0\in A^+_{r,q}(\gamma)$ and
$[b_1,I_\alpha^{\gamma+}]_k$ is bounded from $L^r(\omega_0^r)$ to
$L^q(\omega_0^q)$, whereas $[b_2,I_\alpha^{\gamma+}]_k$ is compact from
$L^r(\omega_0^r)$ to $L^q(\omega_0^q)$ (see Remarks
\ref{rmk:20260729.1619} and \ref{rmk:20260729.1620}).

These examples indicate a fundamental difference from the classical
Euclidean case, where the weighted boundedness and the weighted compactness
of commutators of fractional integral operators with respect to
a single classical Muckenhoupt $A_{r,q}$ weight can indeed characterize the classical
$\mathrm{BMO}$ and $\mathrm{CMO}$ (see \cite{amr(iumj-2020),gwy(cjm-2021)}).
\end{remark}

The implications (ii)$\implies$(iii)$\implies$(iv)$\implies$(v) in Theorem
\ref{thm:20260714.2255} are trivial. To show the implication (v)$\implies$(i) in
Theorem \ref{thm:20260714.2255}, we need to use Corollary \ref{ExtraCC}, give a
quantitative Uchiyama-type characterization of the distance from functions in PBMO to
PCMO (see Lemma \ref{Uchiyama}), establish a local lower estimate for
$[b,I_\alpha^{\gamma+}]_k$ (see Lemma \ref{lem:20260710.1759}), and employ
the following characterization of the distance from a symbol $b\in\mathrm{PBMO}$
to $\mathrm{PCMO}$ in terms of the operator-norm distance from
$[b,I_\alpha^{\gamma+}]_k$ to the space of compact operators from $L^r$ to $L^q$,
which is the third main result of this article. This characterization, combined with
Corollary \ref{ExtraCC}, also provides the proof of (i)$\implies$(ii).

\begin{theorem}\label{distance from PBMO to PCMO}
Let $\gamma,\alpha\in(0,1)$, $1<r<q<\infty$ satisfy $\frac1r-\frac1q=\alpha$,
and $k\in\mathbb{N}$. Then, for any $b\in\mathrm{PBMO}$,
\begin{align*}
\mathrm{dist}_\mathrm{PBMO}(b,\mathrm{PCMO})\sim
\inf\left\{\left\|\left[b,I_\alpha^{\gamma+}\right]_k
-T\right\|_{L^r\to L^q}^\frac{1}{k}:
T\mbox{\ is compact\ from\ }L^r\mbox{\ to\ }L^q\right\},
\end{align*}
where the positive equivalence constants are independent of $b$.
\end{theorem}

Unlike the classical Fr\'echet--Kolmogorov approach used in showing
unweighted compactness of commutators of classical fractional integral
operators (see, for example, \cite{cdw(pa-2009), wy(pams-2018)}), a key
ingredient in the proof of the upper estimate in
Theorem \ref{distance from PBMO to PCMO} is a truncation and finite-rank
approximation argument inspired by Arai and Nakai \cite{an(tjm-2019)}
(see Lemma \ref{lem:20260715.1703}).

We also obtain the following parabolic weighted interpolation of compactness
by establishing two parabolic weighted Fr\'echet--Kolmogorov
criteria (see Propositions \ref{Frechet--Kolmogorov} and \ref{Frechet--Kolmogorov 2}),
which is of independent interest and allows quasi-Banach target spaces.

\begin{theorem}\label{parabolic interpolation of compactness}
Let $\gamma\in(0,1)$, $1\leq r_0,r_1<\infty$, $0<q_0,q_1<\infty$, and
$\omega_0^{q_0},\omega_1^{q_1}\in A_\infty^+(\gamma)$. Let $T$ be a
linear operator defined on $L^{r_0}(\omega_0^{r_0})+L^{r_1}(\omega_1^{r_1})$
and taking values in $\mathscr{M}$. If $T$ is bounded from
$L^{r_0}(\omega_0^{r_0})$ to $L^{q_0}(\omega_0^{q_0})$ and compact from
$L^{r_1}(\omega_1^{r_1})$ to $L^{q_1}(\omega_1^{q_1})$, then, for any
$\theta\in(0,1)$, $T$ is compact from $L^{r_\theta}(\omega_\theta^{r_\theta})$
to $L^{q_\theta}(\omega_\theta^{q_\theta})$, where
\begin{align*}
\frac1{r_\theta}=\frac{1-\theta}{r_0}+\frac\theta{r_1},\ \
\frac1{q_\theta}=\frac{1-\theta}{q_0}+\frac\theta{q_1},\mbox{\ \ and\ \ }
\omega_\theta=\omega_0^{1-\theta}\omega_1^\theta.
\end{align*}
\end{theorem}

The organization of the remainder of this article is as follows. In Section
\ref{section2}, we first establish the unweighted compactness of
mixed commutators of  $I_\alpha^{\gamma+}$ (see Lemma \ref{lem:20260715.1703}),
which is the first of the three tools mentioned above. Then
we present the above-mentioned Uchiyama-type characterization (see Lemma
\ref{Uchiyama}). After that, we use the median method to show a local lower
estimate for $[b,I_\alpha^{\gamma+}]_k$ (see Lemma \ref{lem:20260710.1759}
and Corollary \ref{cor:20260713.2241}). We also construct a counterexample
showing that the boundedness with respect to a single parabolic Muckenhoupt weight
does not characterize PBMO. Finally, we combine the results established above
to prove Theorem \ref{distance from PBMO to PCMO}. In Section \ref{section3},
we first show an interpolation lemma between parabolic weighted Lebesgue spaces
(see Lemma \ref{compact extrapolation lem}). Using this, we then show Theorem
\ref{parabolic extrapolation of compactness}. Combining Theorems
\ref{parabolic extrapolation of compactness} and
\ref{distance from PBMO to PCMO}, we prove Theorem \ref{thm:20260714.2255}.
In Section \ref{section4}, we first establish the parabolic weighted
Fr\'echet--Kolmogorov theorem, which characterizes relatively compact subsets
of parabolic weighted Lebesgue spaces (see Propositions
\ref{Frechet--Kolmogorov} and \ref{Frechet--Kolmogorov 2}). Finally we show
Theorem \ref{parabolic interpolation of compactness}.

We end this introduction by making some notational conventions. Throughout this
article, let $\mathbb{N}:=\{1,2,\ldots\}$ and
$\mathbb{Z}_+:=\mathbb{N}\cup\{0\}$. For any $q\in[1,\infty]$, let $q'$ be the
\emph{conjugate number} of $q$, i.e., $\frac1q+\frac{1}{q'}=1$.
Let $\mathbf{0}$ denote the origin of $\mathbb{R}^{n}$.
For any $A\subset\mathbb{R}^{n+1}$ and $a\in\mathbb{R}$, let
\begin{align*}
A+(\mathbf{0},a):=\{(x,t+a):(x,t)\in A\}\mbox{\ \ and\ \ }
A-(\mathbf{0},a):=\{(x,t-a):(x,t)\in A\}.
\end{align*}
For any measurable set $A\subset\mathbb{R}^{n+1}$, we denote
by $\vert A\vert$ its $(n+1)$-dimensional Lebesgue measure.
For any weight $\omega$ and any measurable set
$A\subset\mathbb{R}^{n+1}$, let
$\omega(A):=\int_A\omega$. The \emph{notation} $f\lesssim g$
means that there exists a positive constant $C$ such that $f\leq Cg$.
If $f\lesssim g$ and $g\lesssim f$, then we write $f\sim
g$. Finally, in all proofs we consistently retain the
notation introduced in the original theorem (or related statement).

\section{Proof of Theorem \ref{distance from PBMO to PCMO}}
\label{section2}

In this section, we prove Theorem \ref{distance from PBMO to PCMO}. To this
end, we first recall the weighted boundedness result for both the operator
$I_\alpha^{\gamma+}$ in \cite[Theorem 5.7]{kyyz(cjm-2025)} and its commutator
in \cite[Theorem 1.9]{ckyyz(jga-2026)}.

\begin{lemma}\label{[I_alpha,b]}
Let $\gamma,\alpha\in(0,1)$. Then the following assertions hold.
\begin{enumerate}
\item[\rm(i)] Let $1<r<q<\infty$ satisfy $\frac1r-\frac1q=\alpha$, and let $\omega$
be a weight. Then $I_\alpha^{\gamma+}$ is bounded from $L^r(\omega^r)$
to $L^q(\omega^q)$ if and only if $\omega\in A_{r,q}^+(\gamma)$.

\item[\rm(ii)] Let $1<r<q<\infty$ satisfy $\frac1r-\frac1q=\alpha$,
$k\in\mathbb{N}$, $b\in\mathrm{PBMO}$, and $\omega\in A_{r,q}^+(\gamma)$. Then
$[b,I_\alpha^{\gamma+}]_k$ is bounded from $L^r(\omega^r)$ to $L^q(\omega^q)$.
Moreover, $\|[b,I_\alpha^{\gamma+}]_k\|_{L^r(\omega^r)\to
L^q(\omega^q)}\lesssim\|b\|_\mathrm{PBMO}^k$, where the implicit positive
constant is independent of $b$.

\item[\rm(iii)] Conversely, if there exist $1<r_0<q_0<\infty$, satisfying
$\frac{1}{r_0}-\frac{1}{q_0}=\alpha$, and $b\in L_\mathrm{loc}^1$ such that
$[b,I_\alpha^{\gamma+}]$ is bounded from $L^{r_0}$ to $L^{q_0}$, then
$b\in\mathrm{PBMO}$.
\end{enumerate}
\end{lemma}

To prove Theorem \ref{distance from PBMO to PCMO}, we also need the following
weighted boundedness estimate for mixed commutators of linear operators,
which follows from an argument similar to the proof of \cite[Theorem
3.9]{ckyyz(jga-2026)}, except that the Cauchy integral formula used therein is
replaced by its multi-parameter version; we omit the details.

\begin{lemma}\label{JGA-Theorem 3.9-Improve}
Let $\gamma\in(0,1)$, $1<r\leq q<\infty$, and
$T:L_{\mathrm{loc}}^1\to\mathscr{M}$ be a linear operator. If, for any
$\omega\in A_{r,q}^+(\gamma)$, $T$ is bounded from
$L^r(\omega^r)$ to $L^q(\omega^q)$ with its operator norm dominated by a
bound only depending on $n,p,\gamma,r,q$, and $[\omega]_{A_{r,q}^+(\gamma)}$,
then, for any  $k\in\mathbb{N}$, $\omega\in A_{r,q}^+(\gamma)$, and
$b_1,\dots,b_k\in\mathrm{PBMO}$, $[b_1,\dots,[b_k,T]]$ is bounded from
$L^r(\omega^r)$ to $L^q(\omega^q)$ and
\begin{align*}
\|[b_1,\dots,[b_k,T]]\|_{L^r(\omega^r)\to L^q(\omega^q)}
\lesssim\prod_{j=1}^k\left\|b_j\right\|_\mathrm{PBMO},
\end{align*}
where the implicit positive constant is independent of $b_1,\dots,b_k$.
\end{lemma}

Then we establish the following compactness result for mixed commutators
of $I_\alpha^{\gamma+}$. It is worth noting that the parabolic domain
$\Omega_{(x,t)}^{\gamma+}$ depends on the base point $(x,t)$. As a result,
unlike the classical Fr\'echet--Kolmogorov approach used for commutators of
classical fractional integral operators (see, for example, \cite{cdw(pa-2009),
wy(pams-2018)}), in the present parabolic setting we employ a truncation and
finite-rank approximation argument inspired by \cite{an(tjm-2019)}. In what
follows, $\varepsilon\to0^+$ means that $\varepsilon\in(0,\infty)$ and
$\varepsilon\to0$.

\begin{lemma}\label{lem:20260715.1703}
Let $\gamma,\alpha\in(0,1)$, $1<r<q<\infty$ satisfy $\frac1r-\frac1q=\alpha$,
$k\in\mathbb{N}$, and $b_1,\ldots,b_k\in\mathrm{PBMO}$. If there exists
$j_0\in\mathbb{N}\cap[1,k]$ such that $b_{j_0}\in\mathrm{PCMO}$, then the
mixed parabolic fractional commutator $[b_1,\ldots,[b_k,I_\alpha^{\gamma+}]]$
is compact from $L^r$ to $L^q$. In particular, for any $b\in\mathrm{PCMO}$,
$[b,I_\alpha^{\gamma+}]_k$ is compact from $L^r$ to $L^q$.
\end{lemma}

\begin{proof}
Since, for any $f\in L^r$ and $(x,t)\in\mathbb{R}^{n+1}$,
\begin{align*}
\left[b_1,\dots,\left[b_k,I_\alpha^{\gamma+}\right]\right]f(x,t)
=\int_{\Omega_{(x,t)}^{\gamma+}}\prod_{j=1}^k\left[b_j(x,t)-b_j(y,s)\right]
\frac{f(y,s)}{[d_p((x,t),(y,s))]^{(n+p)(1-\alpha)}}\,dy\,ds,
\end{align*}
we may assume that $j_0=1$ without loss of generality. We then divide the
proof into the following three steps.

\emph{Step 1}. We first assume that $b_1\in C_\mathrm{c}^\infty$. In this step,
we consider a truncated version $T_\varepsilon$ of the mixed commutator
$[b_1,\dots,[b_k,I_\alpha^{\gamma+}]]$ and show that $T_\varepsilon:L^r\to
L^q$ is compact. To begin with, for any given $\varepsilon\in(0,1)$, define
the truncated kernel $K_\varepsilon$ by setting, for any
$(x,t),(y,s)\in\mathbb{R}^{n+1}$,
\begin{align*}
K_\varepsilon((x,t),(y,s))&:=
\mathbf{1}_{\Omega_{(x,t)}^{\gamma+}}(y,s)
\mathbf{1}_{\{((x,t),(y,s))\in\mathbb{R}^{n+1}\times\mathbb{R}^{n+1}:
\varepsilon<d_p((x,t),(y,s))<\frac{1}{\varepsilon}\}}((x,t),(y,s))\\
&\quad\times\prod_{j=1}^k\left[b_j(x,t)-b_j(y,s)\right]
\frac{1}{[d_p((x,t),(y,s))]^{(n+p)(1-\alpha)}}.
\end{align*}
From H\"older's inequality and the facts that $b_1\in
C_\mathrm{c}^\infty\subset L^\infty$ and
$b_2,\dots,b_k\in\mathrm{PBMO}\hookrightarrow
L_\mathrm{loc}^u$ for any $u\in(1,\infty)$ (which can be deduced directly from
the John--Nirenberg inequality for PBMO functions), it follows that
$K_\varepsilon\in L^q_{(x,t)}L^{r'}_{(y,s)}$. Here, and thereafter,
$L^q_{(x,t)}L^{r'}_{(y,s)}$ denotes the \emph{mixed-norm Lebesgue space}, which is
defined to be the set of all measurable functions $F$ on
$\mathbb{R}^{n+1}\times\mathbb{R}^{n+1}$ such that
\begin{align*}
\|F\|_{L^q_{(x,t)}L^{r'}_{(y,s)}}:=
\left\{\int_{\mathbb{R}^{n+1}}\left[\int_{\mathbb{R}^{n+1}}|F((x,t),(y,s))|
^{r'}\,dy\,ds\right]^\frac{q}{r'}\,dx\,dt\right\}^\frac1q<\infty.
\end{align*}
By the density of product simple functions in $L^q_{(x,t)}L^{r'}_{(y,s)}$
(see, for instance, \cite[Lemma 1.2.19]{Hytonen-book}),
we find that there exists a sequence $\{K_\varepsilon^N\}_{N\in\mathbb{N}}$
of functions in $L^q_{(x,t)}L^{r'}_{(y,s)}$ such that
\begin{align}\label{20260701.2332}
\left\|K_\varepsilon-K_\varepsilon^N\right\|
_{L^q_{(x,t)}L^{r'}_{(y,s)}}\to0\mbox{\ \ as\ \ }N\to\infty.
\end{align}
Moreover, for any $N\in\mathbb{N}$, there exist a sequence
$\{c_j^N\}_{j\in\mathbb{N}\cap[1,N]}$ of real numbers and two sequences
$\{E_j^N\}_{j\in\mathbb{N}\cap[1,N]}$ and $\{F_j^N\}_{j\in\mathbb{N}\cap[1,N]}$
of bounded measurable subsets of $\mathbb{R}^{n+1}$ satisfying, for any
$(x,t),(y,s)\in\mathbb{R}^{n+1}$,
\begin{align*}
K_\varepsilon^N((x,t),(y,s))=\sum_{j=1}^Nc_j^N\mathbf{1}_{E_j^N}(x,t)
\mathbf{1}_{F_j^N}(y,s).
\end{align*}

For any given $\varepsilon\in(0,1)$ and any $N\in\mathbb{N}$, define the
truncated mixed commutators $T_\varepsilon$ and $T_\varepsilon^N$,
respectively, by setting, for any $f\in L^r$ and $(x,t)\in\mathbb{R}^{n+1}$,
\begin{align*}
T_\varepsilon f(x,t)&:=\int_{\mathbb{R}^{n+1}}K_\varepsilon((x,t),(y,s))
f(y,s)\,dy\,ds\\
&\ =\int_{\genfrac{}{}{0pt}{}{\Omega_{(x,t)}^{\gamma+}}
{\varepsilon<d_p((x,t),(y,s))<\frac{1}{\varepsilon}}}
\prod_{j=1}^k\left[b_j(x,t)-b_j(y,s)\right]
\frac{f(y,s)}{[d_p((x,t),(y,s))]^{(n+p)(1-\alpha)}}\,dy\,ds
\end{align*}
and
\begin{align*}
T_\varepsilon^N f(x,t):=\int_{\mathbb{R}^{n+1}}K_\varepsilon^N
((x,t),(y,s))f(y,s)\,dy\,ds=\sum_{j=1}^Nc_j^N\mathbf{1}_{E_j^N}(x,t)
\int_{F_j^N}f(y,s)\,dy\,ds.
\end{align*}
From H\"older's inequality, we deduce that $T_\varepsilon$ is bounded
from $L^r$ to $L^q$ satisfying $\|T_\varepsilon\|_{L^r\to L^q}\leq
\|K_\varepsilon\|_{L^q_{(x,t)}L^{r'}_{(y,s)}}$ and, for any $N\in\mathbb{N}$,
$T_\varepsilon^N$ is a finite-rank operator from $L^r$ to $L^q$.
Applying this, H\"older's inequality and \eqref{20260701.2332}, we obtain
$\|T_\varepsilon-T_\varepsilon^N\|_{L^r\to L^q}\leq
\|K_\varepsilon-K_\varepsilon^N\|_{L^q_{(x,t)}L^{r'}_{(y,s)}}\to0$
as $N\to\infty$, which further implies that
$T_\varepsilon:L^r\to L^q$ is compact.

\emph{Step 2.} We still assume that $b_1\in C_\mathrm{c}^\infty$.
In this step, we prove that
\begin{align}\label{20260720.2302}
\left\|\left[b_1,\dots,\left[b_k,I_\alpha^{\gamma+}\right]\right]
-T_\varepsilon\right\|_{L^r\to L^q}\to0\mbox{\ \ as\ \ }\varepsilon\to0^+,
\end{align}
and hence $[b_1,\dots,[b_k,I_\alpha^{\gamma+}]]:L^r\to L^q$ is compact.
Indeed, applying Lemma \ref{[I_alpha,b]}(i) and Lemma
\ref{JGA-Theorem 3.9-Improve} with $T:=I_\alpha^{\gamma+}$ and $\omega\equiv1$
therein, we find that $[b_1,\dots,[b_k,I_\alpha^{\gamma+}]]$ is bounded from
$L^r$ to $L^q$. Moreover, for any $f\in L^r$ and $(x,t)\in\mathbb{R}^{n+1}$,
\begin{align}\label{20260703.1750}
&\left[b_1,\dots,\left[b_k,I_\alpha^{\gamma+}\right]\right]f(x,t)
-T_\varepsilon f(x,t)\notag\\
&\quad=\left[\int_{\genfrac{}{}{0pt}{}{\Omega_{(x,t)}^{\gamma+}}
{d_p((x,t),(y,s))\leq\varepsilon}}
+\int_{\genfrac{}{}{0pt}{}{\Omega_{(x,t)}^{\gamma+}}
{d_p((x,t),(y,s))\geq\frac{1}{\varepsilon}}}\right]
\prod_{j=1}^k\left[b_j(x,t)-b_j(y,s)\right]
\frac{f(y,s)}{[d_p((x,t),(y,s))]^{(n+p)(1-\alpha)}}\,dy\,ds\notag\\
&\quad=:T^\mathrm{near}_\varepsilon f(x,t)+T^\mathrm{far}_\varepsilon f(x,t).
\end{align}

We first estimate $T^\mathrm{near}_\varepsilon$.
Choose $\delta\in(0,1)$ and $q_\beta\in(q,\infty)$
such that $\beta:=\alpha+\frac{1-\delta}{n+p}\in(0,\frac1r)$ and
$\frac1r-\frac{1}{q_\beta}=\beta$. From the mean value theorem,
it follows that, for any $(x,t),(y,s)\in\mathbb{R}^{n+1}$ with
$d_p((x,t),(y,s))\leq\varepsilon<1$,
\begin{align*}
|b_1(x,t)-b_1(y,s)|&\leq\|\nabla b_1\|_{L^\infty}|x-y|
+\|\partial_tb_1\|_{L^\infty}|t-s|\\
&\lesssim\|\nabla b_1\|_{L^\infty}\|x-y\|_\infty+
\|\partial_tb_1\|_{L^\infty}|t-s|^\frac1p\\
&\leq(\|\nabla b_1\|_{L^\infty}+\|\partial_tb_1\|_{L^\infty})d_p((x,t),(y,s))\\
&\leq(\|\nabla b_1\|_{L^\infty}+\|\partial_tb_1\|_{L^\infty})
\varepsilon^\delta\left[d_p((x,t),(y,s))\right]^{1-\delta}.
\end{align*}
This further implies that, for any $f\in L^r$ and $(x,t)\in\mathbb{R}^{n+1}$,
\begin{align}\label{20260716.1750}
\left|T^\mathrm{near}_\varepsilon f(x,t)\right|&\lesssim
\varepsilon^\delta(\|\nabla b_1\|_{L^\infty}
+\|\partial_tb_1\|_{L^\infty})\notag\\
&\quad\times\int_{\Omega_{(x,t)}^{\gamma+}}
\prod_{j=2}^k\left|b_j(x,t)-b_j(y,s)\right|
\frac{|f(y,s)|}{[d_p((x,t),(y,s))]^{(n+p)(1-\beta)}}\,dy\,ds.
\end{align}
If $k=1$, then any product indexed by $j=2,\dots,k$ is understood to equal
$1$, which, together with Lemma \ref{[I_alpha,b]}(i), further implies that
\begin{align*}
\left\|T^\mathrm{near}_\varepsilon\right\|_{L^r\to L^q}
\lesssim\varepsilon^\delta(\|\nabla b_1\|_{L^\infty}
+\|\partial_tb_1\|_{L^\infty})
\left\|I_\beta^{\gamma+}\right\|_{L^r\to L^{q_\beta}}
\to0\mbox{\ \ as\ \ }\varepsilon\to0^+.
\end{align*}
Without loss of generality, we may suppose that $k\geq2$ and, for any
$j\in\mathbb{N}\cap[2,k]$, $\|b_j\|_\mathrm{PBMO}>0$. Applying the
John--Nirenberg inequality for functions in PBMO, we find that there exists
$\eta\in(0,\infty)$, independent of $b_1,\dots,b_k$, such that, for any
$\sigma_2,\dots,\sigma_k\in\{-1,\,1\}$,
\begin{align*}
e^{\eta\sum_{j=2}^k\frac{\sigma_jb_j}{\|b_j\|_\mathrm{PBMO}}}\in
A_{r,q_\beta}\subset A_{r,q_\beta}^+(\gamma).
\end{align*}
On the other hand, note that, for any $u_2,\dots,u_k\in\mathbb{R}$,
\begin{align*}
\prod_{j=2}^k\left|u_j\right|\lesssim
\prod_{j=2}^k\frac{e^{\eta u_j}+e^{-\eta u_j}}{2}=\frac{1}{2^{k-1}}
\sum_{\sigma_2,\dots,\sigma_k\in\{-1,\,1\}}e^{\eta\sum_{j=2}^k\sigma_ju_j},
\end{align*}
where the implicit positive constant is independent of $u_2,\dots,u_k$.
This further implies that, for any $(x,t),(y,s)\in\mathbb{R}^{n+1}$,
\begin{align*}
\prod_{j=2}^k\left|b_j(x,t)-b_j(y,s)\right|\lesssim
\left(\prod_{j=2}^k\left\|b_j\right\|_\mathrm{PBMO}\right)
\sum_{\sigma_2,\dots,\sigma_k\in\{-1,\,1\}}
e^{\eta\sum_{j=2}^k\frac{\sigma_j[b_j(x,t)-b_j(y,s)]}{\|b_j\|_\mathrm{PBMO}}}.
\end{align*}
Combining this and \eqref{20260716.1750}, we obtain, for any $f\in L^r$ and
$(x,t)\in\mathbb{R}^{n+1}$,
\begin{align}\label{20260716.1810}
\left|T^\mathrm{near}_\varepsilon f(x,t)\right|&\lesssim
\varepsilon^\delta(\|\nabla b_1\|_{L^\infty}+\|\partial_tb_1\|_{L^\infty})
\left(\prod_{j=2}^k\left\|b_j\right\|_\mathrm{PBMO}\right)\notag\\
&\quad\times\sum_{\sigma_2,\dots,\sigma_k\in\{-1,\,1\}}
e^{\eta\sum_{j=2}^k\frac{\sigma_jb_j(x,t)}{\|b_j\|_\mathrm{PBMO}}}
I_\beta^{\gamma+}\left(e^{-\eta\sum_{j=2}^k\frac{\sigma_jb_j}
{\|b_j\|_\mathrm{PBMO}}}|f|\right)(x,t).
\end{align}
Moreover, suppose that $\operatorname{supp}(b_1)\subset R(\mathbf{0},0,L_1)$
with $L_1\in(1,\infty)$. Note that, for any given $\varepsilon\in(0,1)$,
any $(x,t)\in\mathbb{R}^{n+1}\setminus R(\mathbf{0},0,L_1+2)$, and any
$(y,s)\in R(x,t,\varepsilon)$, we have $b_1(x,t)=b_1(y,s)=0$. Thus,
$\operatorname{supp}(T^\mathrm{near}_\varepsilon f)\subset
R(\mathbf{0},0,L_1+2)$. From this, H\"older's inequality,
\eqref{20260716.1810}, the proven conclusion that
$e^{\eta\sum_{j=2}^k\frac{\sigma_jb_j}{\|b_j\|_\mathrm{PBMO}}}\in
A_{r,q_\beta}^+(\gamma)$ for any $\sigma_2,\dots,\sigma_k\in\{-1,\,1\}$,
and Lemma \ref{[I_alpha,b]}(i), we infer that, for any $f\in L^r$,
\begin{align*}
\left\|T^\mathrm{near}_\varepsilon f\right\|_{L^q}
&\leq\left|\operatorname{supp}(T^\mathrm{near}_\varepsilon f)\right|
^{\frac1q-\frac{1}{q_\beta}}
\left\|T^\mathrm{near}_\varepsilon f\right\|_{L^{q_\beta}}\\
&\lesssim\varepsilon^\delta(L_1+2)^{(n+p)(\frac1q-\frac{1}{q_\beta})}
(\|\nabla b_1\|_{L^\infty}+\|\partial_tb_1\|_{L^\infty})
\left(\prod_{j=2}^k\left\|b_j\right\|_\mathrm{PBMO}\right)\\
&\quad\times\sum_{\sigma_2,\dots,\sigma_k\in\{-1,\,1\}}
\left\|e^{\eta\sum_{j=2}^k\frac{\sigma_jb_j(x,t)}{\|b_j\|_\mathrm{PBMO}}}
I_\beta^{\gamma+}\left(e^{-\eta\sum_{j=2}^k\frac{\sigma_jb_j}
{\|b_j\|_\mathrm{PBMO}}}|f|\right)\right\|_{L^{q_\beta}}\\
&\lesssim\varepsilon^\delta(L_1+2)^{(n+p)(\frac1q-\frac{1}{q_\beta})}
(\|\nabla b_1\|_{L^\infty}+\|\partial_tb_1\|_{L^\infty})
\left(\prod_{j=2}^k\left\|b_j\right\|_\mathrm{PBMO}\right)\\
&\quad\times\sum_{\sigma_2,\dots,\sigma_k\in\{-1,\,1\}}
\left\|I_\beta^{\gamma+}\right\|_{L^r(e^{r\eta\sum_{j=2}^k\frac{\sigma_jb_j}
{\|b_j\|_\mathrm{PBMO}}})\to L^{q_\beta}(e^{q_\beta\eta\sum_{j=2}^k
\frac{\sigma_jb_j}{\|b_j\|_\mathrm{PBMO}}})}\|f\|_{L^r}.
\end{align*}
Therefore,
\begin{align}\label{20260720.1642}
\left\|T^\mathrm{near}_\varepsilon\right\|_{L^r\to L^q}\to0
\mbox{\ \ as\ \ }\varepsilon\to0^+.
\end{align}

Now, we estimate $T^\mathrm{far}_\varepsilon$. We restrict
$\varepsilon\in(0,\frac{1}{8L_1})$. For any $(x,t),(y,s)\in\mathbb{R}^{n+1}$,
if $d_p((x,t),(y,s))\geq\frac{1}{\varepsilon}$ and $b_1(x,t)-b_1(y,s)\neq0$,
then either $(x,t)\in R(\mathbf{0},0,2L_1)$ and
$(y,s)\in R(\mathbf{0},0,\frac{1}{4\varepsilon})^\complement$ or $(y,s)\in
R(\mathbf{0},0,2L_1)$ and $(x,t)\in R(\mathbf{0},0,\frac{1}{4\varepsilon})
^\complement$. To proceed, we consider the
following two cases for $(x,t)$ and $(y,s)$.

\emph{Case (1)} $(x,t)\in R(\mathbf{0},0,2L_1)$ and
$(y,s)\in R(\mathbf{0},0,\frac{1}{4\varepsilon})^\complement$. In this case,
there exists $m\in\mathbb{Z}_+$ such that $(y,s)\in
A_m:=R(\mathbf{0},0,\frac{2^{m+1}}{4\varepsilon})\setminus
R(\mathbf{0},0,\frac{2^m}{4\varepsilon})$.
Applying the John--Nirenberg inequality for functions
in PBMO, we obtain, for any $j\in\mathbb{N}\cap[2,k]$,
\begin{align}\label{20260720.2154}
&\left\{\int_{R(\mathbf{0},0,2L_1)}\left[\int_{A_m}
\left|b_j(x,t)-\left(b_j\right)_{R(\mathbf{0},0,2L_1)}\right|
^{r'}\,dy\,ds\right]^\frac{q}{r'}\,dx\,dt\right\}^\frac1q\notag\\
&\quad=|A_m|^\frac{1}{r'}\left[\int_{R(\mathbf{0},0,2L_1)}
\left|b_j(x,t)-\left(b_j\right)_{R(\mathbf{0},0,2L_1)}\right|
^q\,dx\,dt\right]^\frac1q\notag\\
&\quad\lesssim|R(\mathbf{0},0,2L_1)|^\frac1q|A_m|^\frac{1}{r'}
\left\|b_j\right\|_\mathrm{PBMO}
\end{align}
and
\begin{align}\label{20260720.2158}
&\left\{\int_{R(\mathbf{0},0,2L_1)}\left[\int_{A_m}
\left|\left(b_j\right)_{R(\mathbf{0},0,\frac{2^{m+2}}{4\varepsilon})}-b_j(y,s)
\right|^{r'}\,dy\,ds\right]^\frac{q}{r'}\,dx\,dt\right\}^\frac1q\notag\\
&\quad\leq|R(\mathbf{0},0,2L_1)|^\frac1q
\left[\int_{R(\mathbf{0},0,\frac{2^{m+2}}{4\varepsilon})}
\left|b_j(y,s)-\left(b_j\right)_{R(\mathbf{0},0,\frac{2^{m+2}}{4\varepsilon})}
\right|^{r'}\,dy\,ds\right]^\frac{1}{r'}\notag\\
&\quad\lesssim|R(\mathbf{0},0,2L_1)|^\frac1q|A_m|^\frac{1}{r'}
\left\|b_j\right\|_\mathrm{PBMO}.
\end{align}
Moreover, since $\varepsilon<\frac{1}{8L_1}$, it follows that
there exists $N\in\mathbb{N}$ such that
$2^{N-1}L_1<\frac{2^{m+2}}{4\varepsilon}\leq2^NL_1$. Thus,
\begin{align*}
\left|\left(b_j\right)_{R(\mathbf{0},0,2L_1)}-
\left(b_j\right)_{R(\mathbf{0},0,\frac{2^{m+2}}{4\varepsilon})}\right|
&\leq\sum_{i=1}^{N-1}\left|\left(b_j\right)_{R(\mathbf{0},0,2^iL_1)}-
\left(b_j\right)_{R(\mathbf{0},0,2^{i+1}L_1)}\right|\\
&\quad+\left|\left(b_j\right)_{R(\mathbf{0},0,2^NL_1)}-
\left(b_j\right)_{R(\mathbf{0},0,\frac{2^{m+2}}{4\varepsilon})}\right|\\
&\leq\sum_{i=1}^{N-1}
\frac{|R(\mathbf{0},0,2^{i+1}L_1)|}{|R(\mathbf{0},0,2^iL_1)|}
\fint_{R(\mathbf{0},0,2^{i+1}L_1)}\left|b_j-\left(b_j\right)
_{R(\mathbf{0},0,2^{i+1}L_1)}\right|\\
&\quad+\frac{|R(\mathbf{0},0,2^NL_1)|}
{|R(\mathbf{0},0,\frac{2^{m+2}}{4\varepsilon})|}\fint_{R(\mathbf{0},0,2^NL_1)}
\left|b_j-\left(b_j\right)_{R(\mathbf{0},0,2^NL_1)}\right|\\
&\lesssim N\left\|b_j\right\|_\mathrm{PBMO}\sim
\left(m+\log_2\frac{1}{L_1\varepsilon}\right)\left\|b_j\right\|_\mathrm{PBMO}.
\end{align*}
Combining this, \eqref{20260720.2154}, and \eqref{20260720.2158},
we conclude that, for any $j\in\mathbb{N}\cap[2,k]$,
\begin{align*}
&\left\|\left[b_j(x,t)-b_j(y,s)\right]
\mathbf{1}_{R(\mathbf{0},0,2L_1)}(x,t)
\mathbf{1}_{A_m}(y,s)\right\|_{L_{(x,t)}^qL_{(y,s)}^{r'}}\\
&\quad\leq\left\{\int_{R(\mathbf{0},0,2L_1)}\left[\int_{A_m}
\left|b_j(x,t)-\left(b_j\right)_{R(\mathbf{0},0,2L_1)}\right|
^{r'}\,dy\,ds\right]^\frac{q}{r'}\,dx\,dt\right\}^\frac1q\\
&\quad\quad+\left\{\int_{R(\mathbf{0},0,2L_1)}\left[\int_{A_m}
\left|\left(b_j\right)_{R(\mathbf{0},0,2L_1)}-
\left(b_j\right)_{R(\mathbf{0},0,\frac{2^{m+2}}{4\varepsilon})}\right|
^{r'}\,dy\,ds\right]^\frac{q}{r'}\,dx\,dt\right\}^\frac1q\\
&\quad\quad+\left\{\int_{R(\mathbf{0},0,2L_1)}\left[\int_{A_m}
\left|\left(b_j\right)_{R(\mathbf{0},0,\frac{2^{m+2}}{4\varepsilon})}
-b_j(y,s)\right|^{r'}\,dy\,ds\right]^\frac{q}{r'}\,dx\,dt\right\}^\frac1q\\
&\quad\lesssim\left(1+m+\log_2\frac{1}{L_1\varepsilon}\right)
|R(\mathbf{0},0,2L_1)|^\frac1q|A_m|^\frac{1}{r'}
\left\|b_j\right\|_\mathrm{PBMO},
\end{align*}
which, together with H\"older's inequality, further implies that
\begin{align}\label{20260720.2224}
&\left\|\prod_{j=2}^k\left[b_j(x,t)-b_j(y,s)\right]
\mathbf{1}_{R(\mathbf{0},0,2L_1)}(x,t)
\mathbf{1}_{A_m}(y,s)\right\|_{L_{(x,t)}^qL_{(y,s)}^{r'}}\notag\\
&\quad\lesssim\left(1+m+\log_2\frac{1}{L_1\varepsilon}\right)^{k-1}
|R(\mathbf{0},0,2L_1)|^\frac1q|A_m|^\frac{1}{r'}
\prod_{j=2}^k\left\|b_j\right\|_\mathrm{PBMO}\notag\\
&\quad\lesssim L_1^\frac{n+p}{q}\left(\frac{2^m}{\varepsilon}\right)
^\frac{n+p}{r'}\left(1+m+\log_2\frac{1}{L_1\varepsilon}\right)^{k-1}
\prod_{j=2}^k\left\|b_j\right\|_\mathrm{PBMO}.
\end{align}

\emph{Case (2)} $(y,s)\in R(\mathbf{0},0,2L_1)$ and $(x,t)\in
R(\mathbf{0},0,\frac{1}{4\varepsilon})^\complement$. In this case,
there exists $m\in\mathbb{Z}_+$ such that $(x,t)\in
A_m:=R(\mathbf{0},0,\frac{2^{m+1}}{4\varepsilon})\setminus
R(\mathbf{0},0,\frac{2^m}{4\varepsilon})$. From an argument similar to that
used in \eqref{20260720.2224}, we deduce that
\begin{align}\label{20260720.2227}
&\left\|\prod_{j=2}^k\left[b_j(x,t)-b_j(y,s)\right]
\mathbf{1}_{A_m}(x,t)\mathbf{1}_{R(\mathbf{0},0,2L_1)}(y,s)
\right\|_{L_{(x,t)}^qL_{(y,s)}^{r'}}\notag\\
&\quad\lesssim\left(1+m+\log_2\frac{1}{L_1\varepsilon}\right)^{k-1}
|R(\mathbf{0},0,2L_1)|^\frac{1}{r'}|A_m|^\frac1q
\prod_{j=2}^k\left\|b_j\right\|_\mathrm{PBMO}\notag\\
&\quad\lesssim L_1^\frac{n+p}{r'}\left(\frac{2^m}{\varepsilon}\right)
^\frac{n+p}{q}\left(1+m+\log_2\frac{1}{L_1\varepsilon}\right)^{k-1}
\prod_{j=2}^k\left\|b_j\right\|_\mathrm{PBMO}.
\end{align}

Using H\"older's inequality, \eqref{20260720.2224}, \eqref{20260720.2227},
and the assumption that $\frac1r-\frac1q=\alpha$,
we find that, for any $f\in L^r$,
\begin{align*}
\left\|T^\mathrm{far}_\varepsilon f\right\|_{L^q}
&\leq\sum_{m\in\mathbb{Z}_+}\left\|\prod_{j=1}^k\left[b_j(x,t)-b_j(y,s)\right]
\frac{\mathbf{1}_{R(\mathbf{0},0,2L_1)}(x,t)\mathbf{1}_{A_m}(y,s)}
{[d_p((x,t),(y,s))]^{(n+p)(1-\alpha)}}\right\|
_{L_{(x,t)}^qL_{(y,s)}^{r'}}\|f\|_{L^r}\\
&\quad+\sum_{m\in\mathbb{Z}_+}\left\|\prod_{j=1}^k
\left[b_j(x,t)-b_j(y,s)\right]
\frac{\mathbf{1}_{A_m}(x,t)\mathbf{1}_{R(\mathbf{0},0,2L_1)}(y,s)}
{[d_p((x,t),(y,s))]^{(n+p)(1-\alpha)}}\right\|
_{L_{(x,t)}^qL_{(y,s)}^{r'}}\|f\|_{L^r}\\
&\lesssim\sum_{m\in\mathbb{Z}_+}\left[\left(\frac{2^m}{\varepsilon}\right)
^{-(n+p)(1-\alpha)}L_1^\frac{n+p}{q}\left(\frac{2^m}{\varepsilon}\right)
^\frac{n+p}{r'}+\left(\frac{2^m}{\varepsilon}\right)^{-(n+p)(1-\alpha)}
L_1^\frac{n+p}{r'}\left(\frac{2^m}{\varepsilon}\right)^\frac{n+p}{q}\right]\\
&\quad\times\|b_1\|_{L^\infty}\left(1+m+\log_2\frac{1}{L_1\varepsilon}\right)^{k-1}
\prod_{j=2}^k\left\|b_j\right\|_\mathrm{PBMO}\|f\|_{L^r}\\
&\leq\left[L_1^\frac{n+p}{q}\sum_{m\in\mathbb{Z}_+}
\left(\frac{2^m}{\varepsilon}\right)^{-\frac{n+p}{q}}\left(1+m+
\log_2\frac{1}{L_1\varepsilon}\right)^{k-1}\right.\\
&\quad\left.+L_1^\frac{n+p}{r'}\sum_{m\in\mathbb{Z}_+}
\left(\frac{2^m}{\varepsilon}\right)^{-\frac{n+p}{r'}}\left(1+m+
\log_2\frac{1}{L_1\varepsilon}\right)^{k-1}\right]\|b_1\|_{L^\infty}
\prod_{j=2}^k\left\|b_j\right\|_\mathrm{PBMO}\|f\|_{L^r}\\
&\lesssim\left(L_1^\frac{n+p}{q}+L_1^\frac{n+p}{r'}\right)
\left(\varepsilon^\frac{n+p}{q}+\varepsilon^\frac{n+p}{r'}\right)
\left(1+\log_2\frac{1}{L_1\varepsilon}\right)^{k-1}\|b_1\|_{L^\infty}
\prod_{j=2}^k\left\|b_j\right\|_\mathrm{PBMO}\|f\|_{L^r}.
\end{align*}
Thus, $\|T^\mathrm{far}_\varepsilon\|_{L^r\to L^q}\to0$ as $\varepsilon\to0^+$.
Combining this, \eqref{20260703.1750}, and \eqref{20260720.1642}, we conclude
that \eqref{20260720.2302} holds.

\emph{Step 3.} Finally, we assume that $b_1\in\mathrm{PCMO}$ and show
that $[b_1,\dots,[b_k,I_\alpha^{\gamma+}]]:L^r\to L^q$ is compact. By the definition
of PCMO, we find that there exists a sequence $\{b_1^m\}_{m\in\mathbb{N}}$ of
$C_\mathrm{c}^\infty$ functions such that $\|b_1-b_1^m\|_\mathrm{PBMO}\to0$ as
$m\to\infty$. By Step 2, we find that, for any $m\in\mathbb{N}$,
$[b_1^m,\dots,[b_k,I_\alpha^{\gamma+}]]$ is
compact from $L^r$ to $L^q$. Furthermore, from Lemma
\ref{JGA-Theorem 3.9-Improve}, we infer that
$[b_1,\dots,[b_k,I_\alpha^{\gamma+}]]$ is bounded from $L^r$ to $L^q$ and
\begin{align*}
&\left\|\left[b_1,\dots,\left[b_k,I_\alpha^{\gamma+}\right]\right]-
\left[b_1^m,\dots,\left[b_k,I_\alpha^{\gamma+}\right]\right]\right\|
_{L^r\to L^q}\\
&\quad=\left\|\left[b_1-b_1^m,\dots,\left[b_k,I_\alpha^{\gamma+}\right]
\right]\right\|_{L^r\to L^q}\\
&\quad\lesssim\left\|b_1-b_1^m\right\|_\mathrm{PBMO}
\prod_{j=2}^k\left\|b_j\right\|_\mathrm{PBMO}\to0
\end{align*}
as $m\to\infty$. Thus, $[b_1,\dots,[b_k,I_\alpha^{\gamma+}]]:L^r\to L^q$ is compact.
This completes the proof of Lemma \ref{lem:20260715.1703}.
\end{proof}

The following lemma presents a quantitative Uchiyama-type characterization of
the distance from functions in PBMO to PCMO. For the classical Euclidean case
regarding CMO, we refer to \cite[Lemma 3]{u(tmj-1978)}. For convenience and
the sake of self-containment, we provide its detailed proof here. In what
follows, for any $R:=R(x,t,L)\in\mathcal{R}$ with $(x,t)\in\mathbb{R}^{n+1}$
and $L\in(0,\infty)$ and for any $\theta\in(0,\infty)$,
let $\theta R:=R(x,t,\theta L)$.

\begin{lemma}\label{Uchiyama}
For any $b\in\mathrm{PBMO}$,
\begin{align*}
&\mathrm{dist}_\mathrm{PBMO}(b,\mathrm{PCMO})\\
&\quad\sim\max\left\{\lim_{L\to0}\sup_{\genfrac{}{}{0pt}{}
{R\in\mathcal{R}}{l(R)\leq L}}\fint_R|b-b_R|,\,\lim_{L\to\infty}
\sup_{\genfrac{}{}{0pt}{}{R\in\mathcal{R}}{l(R)\geq L}}
\fint_R|b-b_R|,\,\lim_{L\to\infty}\sup_{\genfrac{}{}{0pt}{}
{R\in\mathcal{R}}{R\subset R(\mathbf{0},0,L)^\complement}}
\fint_R|b-b_R|\right\},
\end{align*}
where the positive equivalence constants are independent of $b$.
\end{lemma}

\begin{proof}
Let $b\in\mathrm{PBMO}$ and
\begin{align}\label{20261007.1720}
\mathcal{V}(b):=\max\left\{\lim_{L\to0}\sup_{\genfrac{}{}{0pt}{}
{R\in\mathcal{R}}{l(R)\leq L}}\fint_R|b-b_R|,\,\lim_{L\to\infty}
\sup_{\genfrac{}{}{0pt}{}{R\in\mathcal{R}}{l(R)\geq L}}
\fint_R|b-b_R|,\,\lim_{L\to\infty}\sup_{\genfrac{}{}{0pt}{}
{R\in\mathcal{R}}{R\subset R(\mathbf{0},0,L)^\complement}}
\fint_R|b-b_R|\right\}.
\end{align}
We first show that $\mathcal{V}(b)\leq
\mathrm{dist}_\mathrm{PBMO}(b,\mathrm{PCMO})$. To this end, we claim that,
for any $a\in C_\mathrm{c}^\infty$, $\mathcal{V}(a)=0$. Indeed, let $a\in
C_\mathrm{c}^\infty$. Then
\begin{align*}
\lim_{L\to0}\sup_{\genfrac{}{}{0pt}{}
{R\in\mathcal{R}}{l(R)\leq L}}\fint_R|a-a_R|
&\leq\lim_{L\to0}\sup_{\genfrac{}{}{0pt}{}
{R\in\mathcal{R}}{l(R)\leq L}}\fint_R\fint_R|a(x,t)-a(y,s)|\,dy\,ds\,dx\,dt\\
&\lesssim\lim_{L\to0}\sup_{\genfrac{}{}{0pt}{}{R\in\mathcal{R}}
{l(R)\leq L}}(\|\nabla a\|_{L^\infty}+\|\partial_ta\|_{L^\infty})L
\end{align*}
and
\begin{align*}
\lim_{L\to\infty}\sup_{\genfrac{}{}{0pt}{}{R\in\mathcal{R}}
{l(R)\geq L}}\fint_R|a-a_R|
\leq\lim_{L\to\infty}\sup_{\genfrac{}{}{0pt}{}{R\in\mathcal{R}}
{l(R)\geq L}}\frac{2\|a\|_{L^1}}{|R|}=0.
\end{align*}
In addition, since $\mathrm{supp}(a)$ is bounded, it follows that
\begin{align*}
\lim_{L\to\infty}\sup_{\genfrac{}{}{0pt}{}{R\in\mathcal{R}}
{R\subset R(\mathbf{0},0,L)^\complement}}\fint_R|a-a_R|=0.
\end{align*}
Therefore, $\mathcal{V}(a)=0$, and hence the above claim holds.
Furthermore, for any $a\in C_\mathrm{c}^\infty$ and $R\in\mathcal{R}$,
\begin{align*}
\fint_R|b-b_R|\leq\fint_R|(b-a)-(b-a)_R|+\fint_R|a-a_R|
\leq\|b-a\|_\mathrm{PBMO}+\fint_R|a-a_R|.
\end{align*}
This, together with the above claim that $\mathcal{V}(a)=0$, further implies
that $\mathcal{V}(b)\leq\|b-a\|_\mathrm{PBMO}$. Taking the infimum over all
$a\in C_\mathrm{c}^\infty$, we obtain $\mathcal{V}(b)\leq
\mathrm{dist}_\mathrm{PBMO}(b,\mathrm{PCMO})$.

Next, we show that $\mathrm{dist}_\mathrm{PBMO}(b,\mathrm{PCMO})\lesssim
\mathcal{V}(b)$. Let $\varepsilon\in(0,\infty)$. Then there exist
$0<l_\varepsilon<1<L_\varepsilon<\infty$ such that, for any $R\in\mathcal{R}$,
$\fint_R|b-b_R|<\mathcal{V}(b)+\varepsilon$ if
one of the following three conditions holds.
\begin{enumerate}
\item[\rm(i)] $l(R)\leq l_\varepsilon$.

\item[\rm(ii)] $l(R)\geq L_\varepsilon$.

\item[\rm(iii)] $R\subset R(\mathbf{0},0,L_\varepsilon)^\complement$.
\end{enumerate}
We shall construct $a_\varepsilon\in C_\mathrm{c}^\infty$ such that
$\|b-a_\varepsilon\|_\mathrm{PBMO}\lesssim\mathcal{V}(b)+\varepsilon$,
where the implicit positive constant is independent of $b$ and $\varepsilon$.
Choose $\varphi\in C_\mathrm{c}^\infty$ such that $\varphi$ is nonnegative,
$\mathrm{supp}(\varphi)\subset R(\mathbf{0},0,1)$, and $\|\varphi\|_{L^1}=1$.
Fix $\delta\in(0,\frac{l_\varepsilon}{4})$ and define the parabolic dilation
$\varphi_\delta$ of $\varphi$ by setting, for any $(x,t)\in\mathbb{R}^{n+1}$,
$\varphi_\delta(x,t):=\frac{1}{\delta^{n+p}}\varphi(\frac{x}{\delta},
\frac{t}{\delta^p})$. Also define $b_\delta:=b*\varphi_\delta$. Then $b_\delta$
is infinitely differentiable. We first show that
$\|b-b_\delta\|_\mathrm{PBMO}\lesssim\mathcal{V}(b)+\varepsilon$, where
the implicit positive constant is independent of $b,\varepsilon$, and
$\delta$. Indeed,
from (i), the facts that $|\nabla b_\delta|\lesssim
\frac{\mathcal{V}(b)+\varepsilon}{\delta}$ and
$|\partial_t b_\delta|\lesssim\frac{\mathcal{V}(b)+\varepsilon}{\delta^p}$,
and the mean value theorem, we infer that, for any $R\in\mathcal{R}$
with $l(R)\leq\delta$,
\begin{align}\label{20261011.1905}
\fint_R|(b-b_\delta)-(b-b_\delta)_R|&\leq
\fint_R|b-b_R|+\fint_R|b_\delta-(b_\delta)_R|\notag\\
&\lesssim[\mathcal{V}(b)+\varepsilon]\left\{1+\frac{l(R)}{\delta}
+\left(\frac{l(R)}{\delta}\right)^p\right\}
\lesssim\mathcal{V}(b)+\varepsilon.
\end{align}
Moreover, using the fact that $\|\varphi_\delta\|_{L^1}=1$,
$\mathrm{supp}(\varphi_\delta)\subset
R(\mathbf{0},0,\delta)$, $|\varphi_\delta|\leq
\frac{\|\varphi\|_{L^\infty}}{\delta^{n+p}}$, Tonelli's theorem, and (i), we
obtain, for any $R\in\mathcal{R}$ with $l(R)=\delta$,
\begin{align}\label{20261009.1805}
&\fint_R|(b-b_\delta)-(b-b_\delta)_R|\notag\\
&\quad\leq2\fint_R|b-b_\delta|
\leq2\fint_R|b-b_{3R}|+2\fint_R|b_{3R}-b_\delta|\notag\\
&\quad\leq2\fint_R|b-b_{3R}|+2\fint_R\int_{R(\mathbf{0},0,\delta)}
\varphi_\delta(y,s)|b(x-y,t-s)-b_{3R}|\,dy\,ds\,dx\,dt\notag\\
&\quad\lesssim\fint_R|b-b_{3R}|+\frac{1}{\delta^{n+p}}\fint_R
\int_{R(\mathbf{0},0,\delta)}|b(x-y,t-s)-b_{3R}|\,dy\,ds\,dx\,dt\notag\\
&\quad\lesssim\fint_{3R}|b-b_{3R}|\lesssim\mathcal{V}(b)+\varepsilon.
\end{align}
Furthermore, for any $R:=R(x_R,t_R,l(R))\in\mathcal{R}$ with
$(x_R,t_R)\in\mathbb{R}^{n+1}$ and $l(R)\in[\delta,\infty)$, there exist
$M_R\in\mathbb{N}$ and a sequence $\{P_j\}_{j\in\mathbb{N}\cap[1,M_R]}$ of
pairwise disjoint parabolic rectangles such that $P_j\cap R\neq\emptyset$,
$l(P_j)=\delta$ for any $j\in\mathbb{N}\cap[1,M_R]$, and
$R\subset\bigcup_{j\in\mathbb{N}\cap[1,M_R]}P_j\subset
R(x_R,t_R,l(R)+2\delta)\subset3R$. Combining these and \eqref{20261009.1805},
we conclude that, for any $R\in\mathcal{R}$ with $l(R)\in[\delta,\infty)$,
\begin{align*}
\int_R|(b-b_\delta)-(b-b_\delta)_R|&\leq2\int_R|b-b_\delta|
\leq2\sum_{j=1}^{M_R}\int_{P_j}|b-b_\delta|
\lesssim[\mathcal{V}(b)+\varepsilon]\sum_{j=1}^{M_R}\left|P_j\right|\\
&=[\mathcal{V}(b)+\varepsilon]
\left|\bigcup_{j\in\mathbb{N}\cap[1,M_R]}P_j\right|
\lesssim[\mathcal{V}(b)+\varepsilon]|R|.
\end{align*}
Applying this, \eqref{20261011.1905}, and \eqref{20261009.1805} and by taking the
supremum over all $R\in\mathcal{R}$, we find that
\begin{align*}
\|b-b_\delta\|_\mathrm{PBMO}\lesssim\mathcal{V}(b)+\varepsilon.
\end{align*}

Fix $L\in(16L_\varepsilon,\infty)$ and choose $\widetilde{\psi}\in
C_\mathrm{c}^\infty$ such that $\mathrm{supp}(\widetilde{\psi})\subset
R(\mathbf{0},0,2)$, $\widetilde{\psi}\in[0,1]$, and $\widetilde{\psi}\equiv1$
on $R(\mathbf{0},0,1)$. For any $(x,t)\in\mathbb{R}^{n+1}$, define $\psi(x,t):=
\widetilde{\psi}(\frac{x}{L},\frac{t}{L^p})$,
$a_\varepsilon:=[b_\delta-(b_\delta)_{R(\mathbf{0},0,L)}]\psi$, and
$g_\varepsilon:=[b_\delta-(b_\delta)_{R(\mathbf{0},0,L)}](1-\psi)$. Then
$a_\varepsilon\in C_\mathrm{c}^\infty$. From the proven
conclusion that $\|b-b_\delta\|_\mathrm{PBMO}
\lesssim\mathcal{V}(b)+\varepsilon$ and the fact that constant functions have
zero PBMO-norm, we deduce that
\begin{align}\label{20261011.1655}
\|b-a_\varepsilon\|_\mathrm{PBMO}&\leq\|b-b_\delta\|_\mathrm{PBMO}
+\|b_\delta-a_\varepsilon\|_\mathrm{PBMO}\notag\\
&\lesssim\mathcal{V}(b)+\varepsilon+\|g_\varepsilon\|_\mathrm{PBMO}
+\left\|(b_\delta)_{R(\mathbf{0},0,L)}\right\|_\mathrm{PBMO}\notag\\
&=\mathcal{V}(b)+\varepsilon+\|g_\varepsilon\|_\mathrm{PBMO}.
\end{align}
Therefore, to show that $\|b-a_\varepsilon\|_\mathrm{PBMO}\lesssim
\mathcal{V}(b)+\varepsilon$, it suffices to prove that
$\|g_\varepsilon\|_\mathrm{PBMO}\lesssim\mathcal{V}(b)+\varepsilon$. To
achieve this, for any $R\in\mathcal{R}$, we consider the following two
cases for $l(R)$.

\emph{Case (1)} $l(R)\geq\frac{L}{8}$. In this case, since $L>16L_\varepsilon$,
it follows that $l(R)>L_\varepsilon$. By this, the proven conclusion that
$\|b-b_\delta\|_\mathrm{PBMO}\lesssim\mathcal{V}(b)+\varepsilon$, and (ii),
we obtain
\begin{align*}
\fint_R|b_\delta-(b_\delta)_R|&\leq\fint_R
\left|(b-b_\delta)-(b-b_\delta)_R\right|+\fint_R|b-b_R|\\
&\leq\|b-b_\delta\|_\mathrm{PBMO}+\fint_R|b-b_R|
\lesssim\mathcal{V}(b)+\varepsilon.
\end{align*}
This further implies that
\begin{align}\label{20261011.1605}
\fint_R\left|g_\varepsilon-(g_\varepsilon)_R\right|
&\leq2\fint_R\left|g_\varepsilon-\left[(b_\delta)_R
-(b_\delta)_{R(\mathbf{0},0,L)}\right]\right|\notag\\
&=2\fint_R\left|\left[b_\delta-(b_\delta)_{R(\mathbf{0},0,L)}
-a_\varepsilon\right]-\left[(b_\delta)_R
-(b_\delta)_{R(\mathbf{0},0,L)}\right]\right|\notag\\
&\leq2\fint_R|b_\delta-(b_\delta)_R|+2\fint_R|a_\varepsilon|\notag\\
&\lesssim\mathcal{V}(b)+\varepsilon+\frac{1}{|R|}\int_{R(\mathbf{0},0,2L)}
\left|b_\delta-(b_\delta)_{R(\mathbf{0},0,L)}\right|\notag\\
&\lesssim[\mathcal{V}(b)+\varepsilon]
\left[1+\frac{|R(\mathbf{0},0,2L)|}{|R|}\right]
\lesssim\mathcal{V}(b)+\varepsilon.
\end{align}

\emph{Case (2)} $l(R)<\frac{L}{8}$. In this case, if $R\subset
R(\mathbf{0},0,L)$, then $\fint_R|g_\varepsilon-(g_\varepsilon)_R|=0$, and,
if $R\subset R(\mathbf{0},0,2L)^\complement$, then the fact that
$g_\varepsilon=b_\delta-(b_\delta)_{R(\mathbf{0},0,L)}$ on $R$, the proven
conclusion that $\|b-b_\delta\|_\mathrm{PBMO}\lesssim
\mathcal{V}(b)+\varepsilon$, and (iii) yield that
\begin{align*}
\fint_R\left|g_\varepsilon-(g_\varepsilon)_R\right|&=
\fint_R|b_\delta-(b_\delta)_R|\leq\fint_R
\left|(b-b_\delta)-(b-b_\delta)_R\right|+\fint_R|b-b_R|\\
&\leq\|b-b_\delta\|_\mathrm{PBMO}+\fint_R|b-b_R|
\lesssim\mathcal{V}(b)+\varepsilon.
\end{align*}
It remains to consider the case that $R\nsubset R(\mathbf{0},0,L)$ and
$R\cap R(\mathbf{0},0,2L)\neq\emptyset$. In this case, write
$R:=R(x_R,t_R,l(R))$ with $(x_R,t_R)\in\mathbb{R}^{n+1}$. Then
$d_p((x_R,t_R),(\mathbf{0},0))\in[\frac{7L}{8},\frac{17L}{8}]$ (which follows
directly by contradiction), and hence $R\subset
R(\mathbf{0},0,L_\varepsilon)^\complement$. This, together with an argument
similar to that used in \eqref{20261011.1605} with (ii) therein replaced by
(iii), further implies that $\fint_R\left|g_\varepsilon-(g_\varepsilon)_R\right|
\lesssim\mathcal{V}(b)+\varepsilon$.

Let $J$ be the smallest integer such that $2^Jl(R)\geq L$, i.e.,
$2^{J-1}l(R)<L\leq2^Jl(R)$. For any $j\in\mathbb{Z}_+\cap[0,J]$, define
$R_j:=R(x_R,t_R,2^jl(R))$. We claim that, for any $j\in\mathbb{Z}_+\cap[0,J]$,
\begin{align}\label{20261011.1614}
\fint_{R_j}\left|b_\delta-(b_\delta)_{R_j}\right|
\lesssim\mathcal{V}(b)+\varepsilon,
\end{align}
where the implicit positive constant is independent of $b$ and $j$. Indeed,
let $j\in\mathbb{Z}_+\cap[0,J]$. If $2^jl(R)\leq\frac{L}{2}$, then the proven
conclusion that $d_p((x_R,t_R),(\mathbf{0},0))\geq\frac{7L}{8}$ yields
$R_j\subset R(\mathbf{0},0,L_\varepsilon)^\complement$, and hence
\eqref{20261011.1614} holds. Moreover, if $2^jl(R)>\frac{L}{2}$, then
$l(R_j)>L_\varepsilon$, and hence \eqref{20261011.1614} holds. Thus, the above
claim \eqref{20261011.1614} holds for any $j\in\mathbb{Z}_+\cap[0,J]$.
Combining this and the assumption that $2^{J-1}l(R)<L\leq2^Jl(R)$, we obtain
\begin{align*}
\left|(b_\delta)_R-(b_\delta)_{R(\mathbf{0},0,L)}\right|
&\leq\sum_{j=0}^{J-1}\left|(b_\delta)_{R_j}-(b_\delta)_{R_{j+1}}\right|
+\left|(b_\delta)_{R_J}-(b_\delta)_{R(\mathbf{0},0,L)}\right|\\
&\lesssim(J+1)[\mathcal{V}(b)+\varepsilon]
\sim\left[1+\log_2\frac{L}{l(R)}\right][\mathcal{V}(b)+\varepsilon].
\end{align*}
From this and \eqref{20261011.1614} with $j:=0$ therein, we deduce that
\begin{align*}
\fint_R\left|g_\varepsilon-(g_\varepsilon)_R\right|
&\leq2\fint_R\left|(1-\psi)\left[b_\delta
-(b_\delta)_{R(\mathbf{0},0,L)}\right]-(1-\psi_R)
\left[(b_\delta)_R-(b_\delta)_{R(\mathbf{0},0,L)}\right]\right|\\
&\leq2\fint_R|b_\delta-(b_\delta)_R|
+2\left|(b_\delta)_R-(b_\delta)_{R(\mathbf{0},0,L)}\right|\fint_R|\psi-\psi_R|\\
&\lesssim\left[1+\frac{l(R)}{L}\log_2\frac{L}{l(R)}\right]
[\mathcal{V}(b)+\varepsilon]\lesssim\mathcal{V}(b)+\varepsilon.
\end{align*}

Finally, combining the arguments in Cases (1) and (2), we conclude that
$\|g_\varepsilon\|_\mathrm{PBMO}\lesssim\mathcal{V}(b)+\varepsilon$. This,
together with \eqref{20261011.1655}, further implies that
\begin{align*}
\mathrm{dist}_\mathrm{PBMO}(b,\mathrm{PCMO})
\leq\|b-a_\varepsilon\|_\mathrm{PBMO}\lesssim\mathcal{V}(b)+\varepsilon.
\end{align*}
By letting $\varepsilon\to0^+$, we obtain
$\mathrm{dist}_\mathrm{PBMO}(b,\mathrm{PCMO})\lesssim\mathcal{V}(b)$. This
completes the proof of Lemma \ref{Uchiyama}.
\end{proof}

We also need the following local lower estimate for the commutator
$[b,I_\alpha^{\gamma+}]_k$, which follows from a median argument.

\begin{lemma}\label{lem:20260710.1759}
Let $\gamma,\alpha\in(0,1)$, $1<r<q<\infty$ satisfy $\frac1r-\frac1q=\alpha$,
$k\in\mathbb{N}$, and $b\in L_\mathrm{loc}^1$. Then there exists a positive
constant $C$, depending only on $n,p,\gamma,\alpha,r,q$, and $k$, such that,
for any $R\in\mathcal{R}$, there exist a translation $R^*$ of $R$ and a function
$f^R\in L^r$ such that $|f^R|\leq|R|^{-\frac1r}\mathbf{1}_{R^*}$,
$\|f^R\|_{L^r}\leq1$, and
\begin{align*}
\left\|\left[b,I_\alpha^{\gamma+}\right]_kf^R\right\|_{L^q}
\geq C\left(\fint_R|b-b_R|\right)^k.
\end{align*}
\end{lemma}

\begin{proof}
Without loss of generality, we may assume that $\int_R|b-b_R|>0$.
Let $R:=R(x_0,t_0,L_0)\in\mathcal{R}$ with $(x_0,t_0)\in\mathbb{R}^{n+1}$ and
$L_0\in(0,\infty)$. Fix $\Theta\in(0,\infty)$ with
$\Theta>\max\{2+\gamma2^p,\,\frac{2(1+\gamma)}{1-\gamma}\}$ and define
$R^*:=R(x_0,t_0+\Theta L_0^p,L_0)$. Then, for any $(x,t)\in R$ and $(y,s)\in R^*$,
\begin{align}\label{20261005.1852}
\gamma\|y-x\|_\infty^p\leq\gamma2^pL_0^p<(\Theta-2)L_0^p\leq s-t
<(\Theta+1-\gamma)L_0^p.
\end{align}
In particular, $d_p((x,t),(y,s))\sim L_0$ and $(y,s)\in\Omega_{(x,t)}^{\gamma+}$,
which further implies that $R^*\subset\Omega_{(x,t)}^{\gamma+}$.

Denote by $m_b(R^*)$ a \emph{$\frac12$-median} of $b$ on $R^*$, i.e.,
$m_b(R^*)\in\mathbb{R}$ satisfies
\begin{align*}
\left|R^*\cap\left\{b<m_b\left(R^*\right)\right\}\right|\leq\frac{|R^*|}{2}
\mbox{\ \ and\ \ }
\left|R^*\cap\left\{b>m_b\left(R^*\right)\right\}\right|\leq\frac{|R^*|}{2}.
\end{align*}
Then
\begin{align*}
\left|R^*\cap\left\{b\leq m_b\left(R^*\right)\right\}\right|\geq
\frac{|R^*|}{2}\mbox{\ \ and\ \ }\left|R^*\cap\left\{b\geq
m_b\left(R^*\right)\right\}\right|\geq\frac{|R^*|}{2}.
\end{align*}
Note that
\begin{align*}
\int_R|b-b_R|&\leq2\inf_{c\in\mathbb{R}}\int_R|b-c|
\leq2\int_R\left|b-m_b\left(R^*\right)\right|\\
&\leq2\int_{R\cap\{b\leq m_b(R^*)\}}\left[m_b\left(R^*\right)-b\right]+
2\int_{R\cap\{b\geq m_b(R^*)\}}\left[b-m_b\left(R^*\right)\right].
\end{align*}
Therefore, either
\begin{align}\label{20260713.1727}
\int_{R\cap\{b\leq m_b(R^*)\}}\left[m_b\left(R^*\right)-b\right]\geq
\frac14\int_R|b-b_R|
\end{align}
or
\begin{align}\label{20260713.1728}
\int_{R\cap\{b\geq m_b(R^*)\}}\left[b-m_b\left(R^*\right)\right]\geq
\frac14\int_R|b-b_R|.
\end{align}
Define
\begin{align*}
E:=\begin{cases}
R\cap\left\{b\leq m_b\left(R^*\right)\right\} & \mbox{if\ }
\eqref{20260713.1727}\mbox{\ holds},\\
R\cap\left\{b\geq m_b\left(R^*\right)\right\} & \mbox{if\ }
\eqref{20260713.1728}\mbox{\ holds}
\end{cases}
\mbox{\ \ and\ \ }
F:=\begin{cases}
R^*\cap\left\{b\geq m_b\left(R^*\right)\right\} & \mbox{if\ }
\eqref{20260713.1727}\mbox{\ holds},\\
R^*\cap\left\{b\leq m_b\left(R^*\right)\right\} & \mbox{if\ }
\eqref{20260713.1728}\mbox{\ holds}.
\end{cases}
\end{align*}
Then
\begin{align}\label{20261005.1838}
\int_E\left|b-m_b\left(R^*\right)\right|\geq\frac14\int_R|b-b_R|,
\end{align}
$|F|\geq\frac{|R^*|}{2}$, and, for any $(x,t)\in E$ and $(y,s)\in F$,
$b(x,t)-b(y,s)$ does not change sign and
\begin{align}\label{20261005.1848}
|b(x,t)-b(y,s)|\geq\left|b-m_b\left(R^*\right)\right|.
\end{align}
Since $b$ is finite almost everywhere, it follows that there exists
$N\in\mathbb{N}$ such that $F_N:=F\cap\{|b|\leq N\}$ has measure at least
$\frac{|R^*|}{4}$. Define $f^R:=|R|^{-\frac1r}\boldsymbol{1}_{F_N}$. Then it is
easy to show that $f^R\in L^r$, $|f^R|\leq|R|^{-\frac1r}\mathbf{1}_{R^*}$, and
$\|f^R\|_{L^r}\leq1$. Furthermore, by \eqref{20261005.1852},
\eqref{20261005.1848}, and the fact that $|F_N|\geq\frac{|R^*|}{4}$,
we obtain, for any $(x,t)\in E$,
\begin{align}\label{20261029.1110}
\left|\left[b,I_\alpha^{\gamma+}\right]_kf^R(x,t)\right|
\gtrsim|R|^{\alpha-\frac1r}\left|b(x,t)-m_b\left(R^*\right)\right|^k.
\end{align}
Finally, combining this, H\"older's inequality, \eqref{20261005.1838},
and the assumption that $\alpha=\frac1r-\frac1q$, we conclude that
\begin{align*}
\left\|\left[b,I_\alpha^{\gamma+}\right]_kf^R\right\|_{L^q}
&\geq\left\|\boldsymbol{1}_E\left[b,I_\alpha^{\gamma+}\right]_kf^R\right\|_{L^q}
\gtrsim|R|^{\alpha-\frac1r}
\left[\int_E\left|b-m_b\left(R^*\right)\right|^{kq}\right]^\frac1q\\
&\geq|R|^{\alpha-\frac1r+\frac1q-k}
\left[\int_E\left|b-m_b\left(R^*\right)\right|\right]^k
\gtrsim\left(\fint_R|b-b_R|\right)^k.
\end{align*}
This completes the proof of Lemma \ref{lem:20260710.1759}.
\end{proof}

As a direct corollary of Lemma \ref{lem:20260710.1759}, we obtain the following
result, which is an improvement on Lemma \ref{[I_alpha,b]}(iii) from
first-order commutators to those of arbitrary order and from the unweighted
boundedness to weighted boundedness with some weight belonging to $A_{r,q}^-(\gamma)$.

\begin{corollary}\label{cor:20260713.2241}
Let $\gamma,\alpha\in(0,1)$, and let $b\in L_\mathrm{loc}^1$. If there exist
$1<r_0<q_0<\infty$ satisfying $\frac1{r_0}-\frac1{q_0}=\alpha$, $k_0\in\mathbb{N}$,
and $\omega_0\in A_{r_0,q_0}^-(\gamma)$ such that $[b,I_\alpha^{\gamma+}]_{k_0}$ is
bounded from $L^{r_0}(\omega_0^{r_0})$ to $L^{q_0}(\omega_0^{q_0})$,
then $b\in\mathrm{PBMO}$.
\end{corollary}

\begin{proof}
Let $R:=R(x_0,t_0,L_0)\in\mathcal{R}$ with $(x_0,t_0)\in\mathbb{R}^{n+1}$ and
$L_0\in(0,\infty)$. Let $\Theta$, $R^*$, and $f^R$ be the same as in the proof of
Lemma \ref{lem:20260710.1759} and define
$P:=R(x_0,t_0+\frac{\Theta}{2}L_0^p,(\frac{\Theta}{2}+1)^\frac1pL_0)$.
Then a simple calculation yields that $R\subset P^-(\gamma)$ and $R^*\subset
P^+(\gamma)$ with $|R|\sim|P^-(\gamma)|$. By this, \eqref{20261005.1838}, H\"older's
inequality, \eqref{20261029.1110}, and the assumption that
$\frac1{r_0}-\frac1{q_0}=\alpha$, we obtain
\begin{align*}
\left(\fint_R|b-b_R|\right)^{k_0}&\lesssim
\left[\frac{1}{|R|}\int_E\left|b-m_b\left(R^*\right)\right|\right]^{k_0}
\leq\frac{1}{|R|}\int_E\left|b-m_b\left(R^*\right)\right|^{k_0}\\
&\lesssim|R|^{-1-\alpha}\int_E\left|\left[b,I_\alpha^{\gamma+}\right]_{k_0}
\left(|R|^\frac{1}{r_0}f^R\right)\right|\\
&\leq|R|^{-1-\alpha}\left\|\left[b,I_\alpha^{\gamma+}\right]_{k_0}
\left(|R|^\frac{1}{r_0}f^R\right)\right\|_{L^{q_0}(\omega_0^{q_0})}
\left(\int_R\omega^{-q_0'}\right)^\frac{1}{q_0'}\\
&\leq|R|^{-1-\alpha}\left\|\left[b,I_\alpha^{\gamma+}\right]_{k_0}\right\|
_{L^{r_0}(\omega_0^{r_0})\to L^{q_0}(\omega_0^{q_0})}
\left(\int_{R^*}\omega_0^{r_0}\right)^\frac{1}{r_0}
\left(\int_R\omega_0^{-q_0'}\right)^\frac{1}{q_0'}\\
&\lesssim[\omega_0]_{A_{r_0,q_0}^-(\gamma)}
\left\|\left[b,I_\alpha^{\gamma+}\right]_{k_0}\right\|
_{L^{r_0}(\omega_0^{r_0})\to L^{q_0}(\omega_0^{q_0})}.
\end{align*}
Taking the supremum over all $R\in\mathcal{R}$, we find that $b\in\mathrm{PBMO}$.
This completes the proof of Corollary \ref{cor:20260713.2241}.
\end{proof}

\begin{remark}\label{rmk:20260729.1619}
Let all the notation be as in Corollary \ref{cor:20260713.2241}.
We point out that there exist $b_0\in L_\mathrm{loc}^1\setminus\mathrm{PBMO}$
and a weight $\omega_0$ such that, for any $1<r<q<\infty$ satisfying
$\frac1r-\frac1q=\alpha$ and for any $k\in\mathbb{N}$, $\omega_0\in
A_{r,q}^+(\gamma)$ and $[b_0,I_\alpha^{\gamma+}]_k$ is bounded from
$L^r(\omega_0^r)$ to $L^q(\omega_0^q)$.

Indeed, let $1<r<q<\infty$ satisfy $\frac1r-\frac1q=\alpha$, and let
$k\in\mathbb{N}$. For any $(x,t)\in\mathbb{R}^{n+1}$, define $b_0(x,t):=t$ and
$\omega_0(x,t):=e^t$. Then it is easy to show that $b_0\in
L_\mathrm{loc}^1\setminus\mathrm{PBMO}$ and $\omega_0\in A_{r,q}^+(\gamma)$. On
the other hand, for any $f\in L^r(\omega_0^{r})$ and $(x,t)\in\mathbb{R}^{n+1}$,
\begin{align*}
&\left[b_0,I_\alpha^{\gamma+}\right]_kf(x,t)\omega_0(x,t)\\
&\quad=\int_{\Omega_{(x,t)}^{\gamma+}}\frac{(t-s)^kf(y,s)e^t}
{[d_p((x,t),(y,s))]^{(n+p)(1-\alpha)}}\,dy\,ds\\
&\quad=\int_{\{(y,s)\in\mathbb{R}^{n+1}:s<-\gamma\|y\|_\infty^p\}}
\frac{s^ke^s}{[d_p((y,s),(\mathbf{0},0))]^{(n+p)(1-\alpha)}}
f(x-y,t-s)e^{t-s}\,dy\,ds.
\end{align*}
For any $(y,s)\in\mathbb{R}^{n+1}$, define
\begin{align*}
K_k(y,s):=\frac{s^ke^s}{[d_p((y,s),(\mathbf{0},0))]^{(n+p)(1-\alpha)}}
\mathbf{1}_{\{(y,s)\in\mathbb{R}^{n+1}:s<-\gamma\|y\|_\infty^p\}}(y,s)
\end{align*}
and
\begin{align*}
g(y,s):=e^sf(y,s).
\end{align*}
Then $([b_0,I_\alpha^{\gamma+}]_kf)\omega_0=K_k\ast g$. Moreover, let
$\sigma:=\frac{1}{1-\alpha}$. Then $1+\frac1q=\frac1r+\frac1\sigma$ and
\begin{align*}
\|K_k\|_{L^\sigma}^\sigma
&=\int_{\{(y,s)\in\mathbb{R}^{n+1}:s<-\gamma\|y\|_\infty^p\}}
\frac{|s|^{k\sigma}e^{\sigma s}}{\max\{\|y\|_\infty,\,|s|^\frac1p\}^{n+p}}
\,dy\,ds\\
&\lesssim\int_{-\infty}^0\int_{\{y\in\mathbb{R}^n:\|y\|_\infty
\leq\gamma^{-\frac{1}{p}}|s|^{\frac{1}{p}}\}}
|s|^{k\sigma-1-\frac{n}{p}}e^{\sigma s}\,dy\,ds\\
&\sim\int_{-\infty}^0|s|^{k\sigma-1}e^{\sigma s}\,ds<\infty.
\end{align*}
From this, the fact that $1+\frac1q=\frac1r+\frac1\sigma$, and Young's
inequality, we deduce that, for any $k\in\mathbb{N}$ and $f\in L^r(\omega_0^r)$,
$\|[b_0,I_\alpha^{\gamma+}]_kf\|_{L^q(\omega_0^q)}
=\|K_k\ast g\|_{L^q}\leq\|K_k\|_{L^\sigma}\|g\|_{L^r}\lesssim
\|f\|_{L^r(\omega_0^r)}$, and hence $[b_0,I_\alpha^{\gamma+}]_k$ is bounded
from $L^r(\omega_0^r)$ to $L^q(\omega_0^q)$.

The above example shows a phenomenon that is entirely different from the
classical Euclidean case. More precisely, in the Euclidean setting, if, for
some $k\in\mathbb N$ and some classical Muckenhoupt $A_{r,q}$ weight,
the $k$-th order commutator of the classical fractional integral operator
is bounded on the corresponding weighted space, then the symbol belongs to the
classical space BMO (see, for instance, \cite[Corollary 1.4]{amr(iumj-2020)}
and \cite[Corollary 4.3]{gwy(cjm-2021)}).
\end{remark}

Now, we are ready to prove Theorem \ref{distance from PBMO to PCMO}.

\begin{proof}[Proof of Theorem \ref{distance from PBMO to PCMO}]
Let $b\in\mathrm{PBMO}$. We first show that
\begin{align}\label{2200}
\inf\left\{\left\|\left[b,I_\alpha^{\gamma+}\right]_k-T\right\|
_{L^r\to L^q}^\frac{1}{k}:T\mbox{\ is compact\ from\ }L^r\mbox{\ to\ }L^q\right\}
\lesssim\mathrm{dist}_\mathrm{PBMO}(b,\mathrm{PCMO}).
\end{align}
Let $a\in\mathrm{PCMO}$ and $\widetilde{b}:=b-a$. Then, for any $(x,t),(y,s)\in
\mathbb{R}^{n+1}$,
\begin{align*}
\left[b(x,t)-b(y,s)\right]^k&=\left\{\left[\widetilde{b}(x,t)
-\widetilde{b}(y,s)\right]+\left[a(x,t)-a(y,s)\right]\right\}^k\\
&=\sum_{j=0}^k\dbinom{k}{j}\left[\widetilde{b}(x,t)-\widetilde{b}(y,s)\right]^j
\left[a(x,t)-a(y,s)\right]^{k-j},
\end{align*}
which further implies that, for any $f\in L^r$ and $(x,t)\in\mathbb{R}^{n+1}$,
\begin{align}\label{20260727.2140}
\left[b,I_\alpha^{\gamma+}\right]_kf(x,t)
=\left[\widetilde{b},I_\alpha^{\gamma+}\right]_kf(x,t)
+\sum_{j=0}^{k-1}\dbinom{k}{j}\left[\widetilde{b},
\left[a,I_\alpha^{\gamma+}\right]_{k-j}\right]_jf(x,t).
\end{align}
Furthermore, from Lemma \ref{lem:20260715.1703}, we deduce that
$\sum_{j=0}^{k-1}\binom{k}{j}[\widetilde{b},
[a,I_\alpha^{\gamma+}]_{k-j}]_j$ is compact from $L^r$ to $L^q$,
where we used the notation $[b,T]_0:=T$.
Combining
this, \eqref{20260727.2140}, and Lemma \ref{[I_alpha,b]}(ii), we obtain
\begin{align*}
&\inf\left\{\left\|\left[b,I_\alpha^{\gamma+}\right]_k-T\right\|_{L^r\to L^q}:
T\mbox{\ is compact\ from\ }L^r\mbox{\ to\ }L^q\right\}\\
&\quad\leq\left\|\left[b,I_\alpha^{\gamma+}\right]_k-
\sum_{j=0}^{k-1}\dbinom{k}{j}\left[\widetilde{b},
\left[a,I_\alpha^{\gamma+}\right]_{k-j}\right]_j\right\|_{L^r\to L^q}\\
&\quad=\left\|\left[\widetilde{b},I_\alpha^{\gamma+}\right]_k\right\|
_{L^r\to L^q}\lesssim\left\|b-a\right\|_\mathrm{PBMO}^k.
\end{align*}
Taking the infimum over all $a\in\mathrm{PCMO}$ yields \eqref{2200}.

Next, we show that
\begin{align}\label{20260728.2007}
\mathrm{dist}_\mathrm{PBMO}(b,\mathrm{PCMO})
\lesssim\inf\left\{\left\|\left[b,I_\alpha^{\gamma+}\right]_k
-T\right\|_{L^r\to L^q}^\frac{1}{k}:T\mbox{\ is compact\ from\ }
L^r\mbox{\ to\ }L^q\right\}.
\end{align}
Without loss of generality, we may assume that
$\mathrm{dist}_\mathrm{PBMO}(b,\mathrm{PCMO})\in(0,\infty)$. Let
$\mathcal{V}(b)$ be the same as in \eqref{20261007.1720}.
Applying Lemma \ref{Uchiyama}, we find that $\mathcal{V}(b)\in(0,\infty)$.
Fix $\varepsilon\in(0,\mathcal{V}(b))$. By passing to a subsequence, we may
assume without loss of generality that there exists a sequence
$\{R_j\}_{j\in\mathbb{N}}$ of parabolic rectangles such that, for any
$j\in\mathbb{N}$, $\fint_{R_j}|b-b_{R_j}|\geq\mathcal{V}(b)-\varepsilon$ and
one of the following three conditions holds:
\begin{enumerate}
\item[\rm(i)] $l(R_j)\to0$ as $j\to\infty$,

\item[\rm(ii)] $l(R_j)\to\infty$ as $j\to\infty$,

\item[\rm(iii)] for any $j\in\mathbb{N}$, $R_j\subset
R(\mathbf{0},0,j)^\complement$ and $\sup_{j\in\mathbb{N}}l(R_j)<\infty$.
\end{enumerate}
Moreover, from Lemma \ref{lem:20260710.1759}, we deduce that there
exists a positive constant $C$ such that, for any
$j\in\mathbb{N}$, one can find a translated parabolic rectangle
$R_j^*$ of $R_j$ and a function $f_j\in L^r$ such that
$|f_j|\leq|R_j|^{-\frac1r}\mathbf{1}_{R_j^*}$,
$\|f_j\|_{L^r}\leq1$, and
\begin{align}\label{20261005.2117}
\left\|\left[b,I_\alpha^{\gamma+}\right]_kf_j\right\|_{L^q}
\geq C\left(\fint_{R_j}\left|b-b_{R_j}\right|\right)^k
\geq C[\mathcal{V}(b)-\varepsilon]^k.
\end{align}
We claim that $f_j\rightharpoonup0$ as $j\to\infty$, i.e., for any given $h\in
L^{r'}$, $\int_{\mathbb{R}^{n+1}}f_jh\to0$ as $j\to\infty$. Using this claim,
we find that, for any compact operator $T:L^r\to L^q$, $\|Tf_j\|_{L^q}\to0$
as $j\to\infty$. Combining this and \eqref{20261005.2117},
we conclude that, for any compact operator $T:L^r\to L^q$,
\begin{align*}
\left\|\left[b,I_\alpha^{\gamma+}\right]_k-T\right\|_{L^r\to L^q}
&\geq\limsup_{j\to\infty}\left\|\left[b,I_\alpha^{\gamma+}\right]_kf_j
-Tf_j\right\|_{L^q}\\
&\geq\limsup_{j\to\infty}\left\{\left\|\left[b,I_\alpha^{\gamma+}\right]_kf_j
\right\|_{L^q}-\left\|Tf_j\right\|_{L^q}\right\}
\gtrsim[\mathcal{V}(b)-\varepsilon]^k.
\end{align*}
Letting $\varepsilon\to0$, taking the infimum over all compact operators
$T:L^r\to L^q$, and using Lemma \ref{Uchiyama}, we then obtain
\eqref{20260728.2007}.

Finally, we prove the above claim.
We consider the following three cases for $R_j$.

\emph{Case (1)} $|R_j^*|=|R_j|\sim[l(R_j)]^{n+p}\to0$ as $j\to\infty$.
In this case, from this, the fact that
$|f_j|\leq|R_j|^{-\frac1r}\mathbf{1}_{R_j^*}$ for any $j\in\mathbb{N}$,
H\"older's inequality, and the absolute continuity of $\|\cdot\|_{L^{r'}}$, it
follows that, for any $h\in L^{r'}$,
\begin{align*}
\left|\int_{\mathbb{R}^{n+1}}f_jh\right|\leq\left|R_j\right|^{-\frac1r}
\int_{R_j^*}|h|\leq\left\|h\mathbf{1}_{R_j^*}\right\|_{L^{r'}}\to0
\mbox{\ \ as\ \ }j\to\infty.
\end{align*}

\emph{Case (2)} $|R_j^*|=|R_j|\sim[l(R_j)]^{n+p}\to\infty$ as $j\to\infty$.
In this case, $|R_j|^{-\frac1r}\to0$ as $j\to\infty$. By this and the fact
that $|f_j|\leq|R_j|^{-\frac1r}\mathbf{1}_{R_j^*}$ for any
$j\in\mathbb{N}$, we obtain, for any $h\in C_\mathrm{c}^\infty$,
\begin{align*}
\left|\int_{\mathbb{R}^{n+1}}f_jh\right|\leq\left|R_j\right|^{-\frac1r}
\int_{R_j^*}|h|\leq\left|R_j\right|^{-\frac1r}\|h\|_{L^1}\to0
\mbox{\ \ as\ \ }j\to\infty.
\end{align*}
This, together with the facts that $C_{\mathrm{c}}^\infty$ is dense in $L^{r'}$ and
that $\{f_j\}_{j\in\mathbb{N}}$ is bounded in $L^r$, further implies that,
for any $h\in L^{r'}$, $\int_{\mathbb{R}^{n+1}}f_jh\to0$ as $j\to\infty$.

\emph{Case (3)} For any $j\in\mathbb{N}$, $R_j\subset
R(\mathbf{0},0,j)^\complement$ and $\sup_{j\in\mathbb{N}}l(R_j)<\infty$.
In this case, we observe that $R_j^*$ tends to infinity as $j\to\infty$.
Therefore, for any $h\in C_\mathrm{c}^\infty$, there exists $j_h\in\mathbb{N}$
such that, for any $j\in\mathbb{N}\cap[j_h,\infty)$,
$\operatorname{supp}(h)\cap R_j^*=\emptyset$. This, together with the fact
that $|f_j|\leq|R_j|^{-\frac1r}\mathbf{1}_{R_j^*}$ for any
$j\in\mathbb{N}$, further implies that, for any
$j\in\mathbb{N}\cap[j_h,\infty)$, $\int_{\mathbb{R}^{n+1}}f_jh=0$, and hence
$\int_{\mathbb{R}^{n+1}}f_jh\to0$ as $j\to\infty$. From this and the facts
that $C_{\mathrm{c}}^\infty$ is dense in $L^{r'}$ and that $\{f_j\}_{j\in\mathbb{N}}$
is bounded in $L^r$, we deduce that, for any $h\in L^{r'}$,
$\int_{\mathbb{R}^{n+1}}f_jh\to0$ as $j\to\infty$.

Combining Cases (1), (2), and (3), we conclude that the above claim
holds, which completes the proof of Theorem \ref{distance from PBMO to PCMO}.
\end{proof}

\section{Proofs of Theorems \ref{parabolic extrapolation of compactness} and \ref{thm:20260714.2255}}
\label{section3}

In this section, we first prove a key interpolation lemma between parabolic
weighted Lebesgue spaces. Then we show Theorems
\ref{parabolic extrapolation of compactness} and \ref{thm:20260714.2255}.
To begin with, we recall the parabolic extrapolation theorem for linear
operators, which was proved in our previous work
\cite[Theorem 1.4]{ckyyz(jga-2026)}.

\begin{lemma}\label{parabolic extrapolation}
Let $\gamma\in(0,1)$ and $1\leq r_0\leq q_0<\infty$. Let $\mathcal{F}$ be the
set of all pairs $(f,g)$ of measurable functions on $\mathbb{R}^{n+1}$ such
that, for any $\omega_0\in A_{r_0,q_0}^+(\gamma)$,
\begin{align}\label{20241204.1150}
\|f\|_{L^{q_0}(\omega_0^{q_0})}\leq C_0\|g\|_{L^{r_0}(\omega_0^{r_0})}
\end{align}
holds provided that the left-hand side of \eqref{20241204.1150} is finite,
where the positive constant $C_0$ in \eqref{20241204.1150} depends only on
$\mathcal{F},n,p,\gamma,r_0,q_0$, and $[\omega_0]_{A_{r_0,q_0}^+(\gamma)}$.
Then, for any $1<r\leq q<\infty$ satisfying $\frac1r-\frac1q
=\frac{1}{r_0}-\frac{1}{q_0}$ and for any $\omega\in A_{r,q}^+(\gamma)$, there
exists a positive constant $C$, depending only on
$\mathcal{F},n,p,\gamma,r_0,q_0,r,q$, and $[\omega]_{A_{r,q}^+(\gamma)}$, such
that, for any $(f,g)\in\mathcal{F}$,
\begin{align}\label{20241204.1152}
\|f\|_{L^q(\omega^q)}\leq C\|g\|_{L^r(\omega^r)}
\end{align}
holds provided that the left-hand side of \eqref{20241204.1152} is finite.
\end{lemma}

Now, we show the above-mentioned interpolation lemma between parabolic
weighted Lebesgue spaces. In what follows, for any $1\leq q_0,q_1<\infty$,
weights $\omega_0,\omega_1$, and $\theta\in(0,1)$, we use the \emph{notation
$[L^{q_0}(\omega_0),L^{q_1}(\omega_1)]_\theta$} to denote the \emph{complex
interpolation space} with respect to $(L^{q_0}(\omega_0),L^{q_1}(\omega_1))$
(see, for instance, \cite[p.\,114, 3]{bl76-book} for its definition). According to
the Stein--Weiss interpolation theorem (see, for instance,
\cite[Corollary 5.5.2]{bl76-book}), we find that
\begin{align*}
\left[L^{q_0}(\omega_0),L^{q_1}(\omega_1)\right]_\theta=L^{q_\theta}(\omega_\theta),
\end{align*}
where $q_\theta\in[1,\infty)$ satisfies
$\frac{1}{q_\theta}=\frac{1-\theta}{q_0}+\frac{\theta}{q_1}$ and where
$\omega_\theta$ is a weight satisfying
$\omega_\theta^\frac{1}{q_\theta}=\omega_0^\frac{1-\theta}{q_0}
\omega_1^\frac{\theta}{q_1}$.

\begin{lemma}\label{compact extrapolation lem}
Let $\gamma\in(0,1)$, $1<r\leq q<\infty$ and $1<r_1\leq q_1<\infty$
satisfy $\frac1r-\frac1q=\frac{1}{r_1}-\frac{1}{q_1}$,
$\omega\in A_{r,q}^+(\gamma)$, and $\omega_1\in A_{r_1,q_1}^-(\gamma)$. Then
there exist $\theta\in(0,1)$, $1<r_0\leq q_0<\infty$, and $\omega_0\in
A_{r_0,q_0}^+(\gamma)$ such that $\frac{1}{r_0}-\frac{1}{q_0}=\frac1r-\frac1q$,
\begin{align*}
L^r(\omega^r)=\left[L^{r_0}\left(\omega_0^{r_0}\right),
L^{r_1}\left(\omega_1^{r_1}\right)\right]_\theta,
\end{align*}
and
\begin{align*}
L^q(\omega^q)=
\left[L^{q_0}\left(\omega_0^{q_0}\right),
L^{q_1}\left(\omega_1^{q_1}\right)\right]_\theta.
\end{align*}
\end{lemma}

\begin{proof}
Since the parabolic Muckenhoupt class is independent of the choice of the time
lag $\gamma$ (see, for instance, \cite[Theorem 3.1]{km(am-2024)}), we may
assume that $\gamma\in(\frac12,1)$ without loss of generality.
Let $\theta\in(0,1)$ be such that $1<r_0(\theta)\leq q_0(\theta)<\infty$, where
\begin{align*}
r_0(\theta):=\frac{1-\theta}{\frac1r-\frac{\theta}{r_1}}
\mbox{\ \ and\ \ }
q_0(\theta):=\frac{1-\theta}{\frac1q-\frac{\theta}{q_1}},
\end{align*}
and let $\omega_0(\theta)$ be a nonnegative function on
$\mathbb{R}^{n+1}$ satisfying $\omega=[\omega_0(\theta)]^{1-\theta}
\omega_1^\theta$. According to the Stein--Weiss interpolation
theorem (see, for instance, \cite[Corollary 5.5.2]{bl76-book}),
it suffices to show that there exists $\theta\in(0,1)$ such that
$$
\omega_0(\theta)=\omega^{\frac{1}{1-\theta}}
\omega_1^{-\frac{\theta}{1-\theta}}\in
A_{r_0(\theta),q_0(\theta)}^+(\gamma).
$$

For any $\theta\in(0,1)$, define $\varepsilon(\theta):=1+\frac{\theta q}{r_1'}$,
$\delta(\theta):=1+\frac{\theta r'}{q_1}$, $f(\theta):=
\frac{\varepsilon(\theta)q_0(\theta)}{(1-\theta)q}$, $f_1(\theta):=
\frac{\theta\varepsilon'(\theta)q_0(\theta)}{(1-\theta)r_1'}$,
$g(\theta):=\frac{\delta(\theta)r_0'(\theta)}{(1-\theta)r'}$,
and $g_1(\theta):=\frac{\theta\delta'(\theta)r_0'(\theta)}
{(1-\theta)q_1}$. Then
$$f(\theta)=f_1(\theta)=
\frac{q_0(\theta)(r_1'+\theta q)}{(1-\theta)qr_1'}
\downarrow1
$$
and
$$
g(\theta)=g_1(\theta)=
\frac{(q_1+\theta r')r_0'(\theta)}{(1-\theta)q_1r'}\downarrow1
$$
as $\theta\to0^+$. By the facts that $\omega^q\in
A_{1+\frac{q}{r'}}^+(\gamma)$, $\omega^{-r'}\in A_{1+\frac{r'}{q}}^-(\gamma)$,
$\omega_1^{q_1}\in A_{1+\frac{q_1}{r_1'}}^-(\gamma)$, and $\omega_1^{-r_1'}\in
A_{1+\frac{r_1'}{q_1}}^+(\gamma)$, and the parabolic reverse H\"older property
of parabolic Muckenhoupt weights (see \cite[Corollary 5.3]{km(am-2024)}), we
find that there exists $\tau_0\in(1,\infty)$ such that,
for any $\tau\in(1,\tau_0]$ and $R:=R(x,t,L)\in\mathcal{R}$
with $(x,t)\in\mathbb{R}^{n+1}$ and $L\in(0,\infty)$,
\begin{align}\label{20250114.0946}
&\left[\fint_{R^-(\gamma)}\omega^{\tau q}\right]
^\frac{1}{\tau}\lesssim\fint_{R^-(\gamma)
+(\mathbf{0},(1-\gamma)L^p)}\omega^q,\ \
\left[\fint_{R^+(\gamma)}\omega^{-\tau r'}\right]
^\frac{1}{\tau}\lesssim\fint_{R^+(\gamma)
-(\mathbf{0},(1-\gamma)L^p)}\omega^{-r'},\notag\\
&\left[\fint_{R^+(\gamma)}\omega_1^{\tau q_1}\right]
^\frac{1}{\tau}\lesssim\fint_{R^+(\gamma)
-(\mathbf{0},(1-\gamma)L^p)}\omega_1^{q_1},\mbox{\ \ and\ \ }
\left[\fint_{R^-(\gamma)}\omega_1^{-\tau r_1'}\right]
^\frac{1}{\tau}\lesssim\fint_{R^-(\gamma)
+(\mathbf{0},(1-\gamma)L^p)}\omega_1^{-r_1'},
\end{align}
where the implicit positive constants are all independent of
$R$. Note that there exists $\theta_0\in(0,1)$ such that,
for any $\theta\in(0,\theta_0)$, $1<f(\theta),f_1(\theta),
g(\theta),g_1(\theta)\leq\tau_0$. Let $\theta:=\frac{\theta_0}{2}$,
$r_0:=r_0(\frac{\theta_0}{2})$, $q_0:=q_0(\frac{\theta_0}{2})$,
$\omega_0:=\omega_0(\frac{\theta_0}{2})$,
$\varepsilon:=\varepsilon(\frac{\theta_0}{2})$, and
$\delta:=\delta(\frac{\theta_0}{2})$. From
H\"older's inequality and \eqref{20250114.0946}, we deduce that
\begin{align}\label{20240603.2205}
\left[\fint_{R^-(\gamma)}\omega_0^{q_0}\right]^\frac{1}{q_0}
&=\left[\fint_{R^-(\gamma)}\omega
^{\frac{q_0}{1-\theta}}\omega_1^{-\frac{\theta
q_0}{1-\theta}}\right]^\frac{1}{q_0}
\leq\left[\fint_{R^-(\gamma)}\omega^{\frac{\varepsilon
q_0}{1-\theta}}\right]^\frac{1}{\varepsilon q_0}
\left[\fint_{R^-(\gamma)}\omega_1^{-\frac{\theta\varepsilon'q_0}
{1-\theta}}\right]^\frac{1}{\varepsilon'q_0}\notag\\
&=\left[\fint_{R^-(\gamma)}\omega^{f(\frac{\theta_0}{2})q}
\right]^\frac{1}{\varepsilon q_0}\left[\fint_{R^-(\gamma)}
\omega_1^{-f_1(\frac{\theta_0}{2})r_1'}\right]
^\frac{1}{\varepsilon'q_0}\notag\\
&\lesssim\left[\fint_{R^-(\gamma)+(\mathbf{0},(1-\gamma)L^p)}
\omega^q\right]^\frac{f(\frac{\theta_0}{2})}{\varepsilon q_0}
\left[\fint_{R^-(\gamma)+(\mathbf{0},(1-\gamma)L^p)}
\omega_1^{-r_1'}\right]^\frac{f_1(\frac{\theta_0}{2})}
{\varepsilon'q_0}\notag\\
&=\left[\fint_{R^-(\gamma)+(\mathbf{0},(1-\gamma)L^p)}\omega^q
\right]^\frac{1}{q(1-\theta)}\left[\fint_{R^-(\gamma)
+(\mathbf{0},(1-\gamma)L^p)}\omega_1^{-r_1'} \right]
^\frac{\theta}{(1-\theta)r_1'}
\end{align}
and
\begin{align*}
\left[\fint_{R^+(\gamma)}\omega_0^{-r_0'}\right]^\frac{1}{r_0'}
&=\left[\fint_{R^+(\gamma)}\omega^{-\frac{r_0'}{1-\theta}}
\omega_1^\frac{\theta r_0'}{1-\theta}\right]^\frac{1}{r_0'}
\leq\left[\fint_{R^+(\gamma)}\omega^{-\frac{\delta r_0'}
{1-\theta}}\right]^\frac{1}{\delta r_0'}
\left[\fint_{R^+(\gamma)}\omega_1^\frac{\theta\delta'r_0'}
{1-\theta}\right]^\frac{1}{\delta'r_0'}\\
&=\left[\fint_{R^+(\gamma)}\omega^{-g(\frac{\theta_0}{2})r'}
\right]^\frac{1}{\delta r_0'}
\left[\fint_{R^+(\gamma)}\omega_1^{g_1(\frac{\theta_0}{2})q_1}
\right]^\frac{1}{\delta'r_0'}\\
&\lesssim\left[\fint_{R^+(\gamma)-(\mathbf{0},(1-\gamma)L^p)}
\omega^{-r'}\right]^\frac{g(\frac{\theta_0}{2})}{\delta r_0'}
\left[\fint_{R^+(\gamma)-(\mathbf{0},(1-\gamma)L^p)}
\omega_1^{q_1}\right]^\frac{g_1(\frac{\theta_0}{2})}{\delta'r_0'}\\
&=\left[\fint_{R^+(\gamma)-(\mathbf{0},(1-\gamma)L^p)}
\omega^{-r'}\right]^\frac{1}{(1-\theta)r'}\left[
\fint_{R^+(\gamma)-(\mathbf{0},(1-\gamma)L^p)}\omega_1^{q_1}
\right]^\frac{\theta}{(1-\theta)q_1}.
\end{align*}
Combining this, \eqref{20240603.2205},
$R^-(\gamma)+(\mathbf{0},(1-\gamma)L^p)\subset R^-(2\gamma-1)$ having equivalent
measures, $R^+(\gamma)-(\mathbf{0},(1-\gamma)L^p)\subset R^+(2\gamma-1)$ having
equivalent measures, and \cite[Theorem 3.1]{km(am-2024)}, we obtain
\begin{align*}
&\left[\fint_{R^-(\gamma)}\omega_0^{q_0}\right]^\frac{1}{q_0}
\left[\fint_{R^+(\gamma)}\omega_0^{-r_0'}\right]^\frac{1}{r_0'}\\
&\quad\lesssim
\left[\fint_{R^-(\gamma)+(\mathbf{0},(1-\gamma)L^p)}\omega^q
\right]^\frac{1}{q(1-\theta)}\left[\fint_{R^-(\gamma)
+(\mathbf{0},(1-\gamma)L^p)}\omega_1^{-r_1'} \right]
^\frac{\theta}{(1-\theta)r_1'}\\
&\quad\quad\times
\left[\fint_{R^+(\gamma)-(\mathbf{0},(1-\gamma)L^p)}
\omega^{-r'}\right]^\frac{1}{(1-\theta)r'}\left[
\fint_{R^+(\gamma)-(\mathbf{0},(1-\gamma)L^p)}\omega_1^{q_1}
\right]^\frac{\theta}{(1-\theta)q_1}\\
&\quad\lesssim\left[\fint_{R^-(2\gamma-1)}\omega^q\right]
^\frac{1}{q(1-\theta)}\left[\fint_{R^+(2\gamma-1)}\omega^{-r'}
\right]^\frac{1}{(1-\theta)r'}\\
&\quad\quad\times\left[\fint_{R^-(2\gamma-1)}\omega_1^{-r_1'}
\right]^\frac{\theta}{(1-\theta)r_1'}\left[\fint_{R^+(2\gamma-1)}
\omega_1^{q_1}\right]^\frac{\theta}{(1-\theta)q_1}\\
&\quad\leq[\omega]_{A_{r,q}^+(2\gamma-1)}^\frac{1}{1-\theta}
[\omega_1]_{A_{r_1,q_1}^-(2\gamma-1)}^\frac{\theta}{(1-\theta)}
<\infty.
\end{align*}
Finally, taking the supremum over all $R\in\mathcal{R}$, we conclude that
$\omega_0\in A_{r_0,q_0}^{+}(\gamma)$. This completes the proof
of Lemma \ref{compact extrapolation lem}.
\end{proof}

We are ready to prove Theorem \ref{parabolic extrapolation of compactness} by
drawing some ideas from \cite[Theorem 2.3]{hl(rmi-2023)}.

\begin{proof}[Proof of Theorem \ref{parabolic extrapolation of compactness}]
Let $1<r\leq q<\infty$ satisfy $\frac1r-\frac1q=\frac1{r_0}-\frac1{q_0}$, and let
$\omega\in A_{r,q}^+(\gamma)$. Lemma \ref{compact extrapolation lem}
yields that there exist $\theta\in(0,1)$, $1<r_\theta\leq
q_\theta<\infty$ with $\frac{1}{r_\theta}-\frac{1}{q_\theta}=
\frac{1}{r}-\frac{1}{q}$, and $\omega_\theta\in
A_{r_\theta,q_\theta}^+(\gamma)$ such that
\begin{align}\label{20261011.1847}
L^r(\omega^r)=\left[L^{r_\theta}\left(\omega_\theta^{r_\theta}
\right),L^{r_1}\left(\omega_1^{r_1}\right)\right]_\theta
\mbox{\ \ and\ \ }L^q(\omega^q)=
\left[L^{q_\theta}\left(\omega_\theta^{q_\theta}\right),
L^{q_1}\left(\omega_1^{q_1}\right)\right]_\theta.
\end{align}
In addition, for any bounded function $h$ on $\mathbb{R}^{n+1}$ with compact
support, define
\begin{align*}
h_N:=\min\{|Th|,\,N\}\boldsymbol{1}_{R(\mathbf{0},0,N)}.
\end{align*}
Then
\begin{align*}
\mathcal{F}:=\left\{(h_N,|h|):h\mbox{\ is \ bounded\ with\ compact\
support\ and\ }N\in\mathbb{N}\right\}
\end{align*}
satisfies the hypothesis of Lemma \ref{parabolic extrapolation}. Applying
Lemma \ref{parabolic extrapolation} and Fatou's lemma as $N\to\infty$,
we obtain, for any bounded function $h$ on
$\mathbb{R}^{n+1}$ with compact support,
$\|Th\|_{L^{q_\theta}(\omega_\theta^{q_\theta})}
\lesssim\|h\|_{L^{r_\theta}(\omega_\theta^{r_\theta})}$. This, together with
a standard density argument, further implies that
$T:L^{r_\theta}(\omega_\theta^{r_\theta})\to
L^{q_\theta}(\omega_\theta^{q_\theta})$ is bounded. Combining this,
\eqref{20261011.1847}, the assumption that $T$ is compact from
$L^{r_1}(\omega_1^{r_1})$ to $L^{q_1}(\omega_1^{q_1})$, and the Cwikel--Kalton
theorem (see, for instance, \cite[Theorem 3.1]{hl(rmi-2023)}), we conclude that
$T$ is compact from $L^r(\omega^r)$ to $L^q(\omega^q)$, which completes the
proof of Theorem \ref{parabolic extrapolation of compactness}.
\end{proof}

We are also ready to show Theorem \ref{thm:20260714.2255}, which is a direct
consequence of Lemma \ref{[I_alpha,b]}(i), Theorem \ref{distance from PBMO to PCMO},
and Corollaries \ref{ExtraCC} and \ref{cor:20260713.2241}.

\begin{proof}[Proof of Theorem \ref{thm:20260714.2255}]
Combining Lemma \ref{[I_alpha,b]}(i), Theorem \ref{distance from PBMO to PCMO}, and
Corollary \ref{ExtraCC}, we find that (i) implies (ii). The implications
(ii)$\implies$(iii)$\implies$(iv)$\implies$(v) are trivial. We then prove
(v)$\implies$(i). From Corollary \ref{cor:20260713.2241}, we deduce that
$b\in\mathrm{PBMO}$. Moreover, using Lemma \ref{[I_alpha,b]}(i), Corollary
\ref{ExtraCC}, and (v), we conclude that $[b,I_\alpha^{\gamma+}]_{k_0}$ is compact
from $L^{r_0}$ to $L^{q_0}$. Finally, this, together with $b\in\mathrm{PBMO}$ and
Theorem \ref{distance from PBMO to PCMO}, further implies that $b\in\mathrm{PCMO}$.
This completes the proof of Theorem \ref{thm:20260714.2255}.
\end{proof}

\begin{remark}\label{rmk:20260729.1620}
Let all the notation be as in Theorem \ref{thm:20260714.2255}.
We now claim that there exist $b_0\in\mathrm{PBMO}\setminus\mathrm{PCMO}$
and a weight $\omega_0$ such that, for any $1<r<q<\infty$
satisfying $\frac1r-\frac1q=\alpha$ and for any $k\in\mathbb{N}$, $\omega_0\in
A_{r,q}^+(\gamma)\setminus A_{r,q}^-(\gamma)$ and $[b_0,I_\alpha^{\gamma+}]_k$ is
compact from $L^r(\omega_0^r)$ to $L^q(\omega_0^q)$.

To prove this, for any $(x,t)\in\mathbb{R}^{n+1}$, define
\begin{align*}
b_0(x,t):=\ln\left(e+d_p((x,t),(\mathbf{0},0))\right)
\mbox{\ \ and\ \ }\omega_0(x,t):=e^t.
\end{align*}
Then it is easy to prove that $\omega_0\in A_{r,q}^+(\gamma)$ and
$b_0\in\mathrm{PBMO}$. This, together with Lemma \ref{[I_alpha,b]}(ii),
further implies that $[b_0,I_\alpha^{\gamma+}]_k$ is bounded from
$L^r(\omega_0^r)$ to $L^q(\omega_0^q)$. However, $b_0\notin\mathrm{PCMO}$.
To see this, we first show that
\begin{align}\label{20260714.1631}
\liminf_{L\to\infty}\fint_{R(\mathbf{0},0,L)}
\left|b_0-(b_0)_{R(\mathbf{0},0,L)}\right|>0.
\end{align}
Observe that, for any $(x,t)\in\mathbb{R}^{n+1}$ and $L\in(0,\infty)$,
\begin{align*}
b_0\left(Lx,L^pt\right)=\ln L+\ln\left(d_p((x,t),(\mathbf{0},0))
+\frac{e}{L}\right),
\end{align*}
and hence
\begin{align*}
(b_0)_{R(\mathbf{0},0,L)}=\ln L+\fint_{R(\mathbf{0},0,1)}
\ln\left(d_p((y,s),(\mathbf{0},0))+\frac{e}{L}\right)\,dy\,ds.
\end{align*}
Therefore, for any given $L\in(0,\infty)$,
\begin{align*}
&\fint_{R(\mathbf{0},0,L)}\left|b_0(x,t
)-(b_0)_{R(\mathbf{0},0,L)}\right|\,dx\,dt\\
&\quad=\fint_{R(\mathbf{0},0,1)}\left|\ln\left(d_p((x,t),(\mathbf{0},0))
+\frac{e}{L}\right)-\fint_{R(\mathbf{0},0,1)}
\ln\left(d_p((y,s),(\mathbf{0},0))+\frac{e}{L}\right)\,dy\,ds\right|\,dx\,dt.
\end{align*}
This, together with the Lebesgue dominated convergence theorem,
further implies that
\begin{align*}
&\lim_{L\to\infty}\fint_{R(\mathbf{0},0,L)}
\left|b_0-(b_0)_{R(\mathbf{0},0,L)}\right|\\
&\quad=\fint_{R(\mathbf{0},0,1)}\left|\ln d_p((x,t),(\mathbf{0},0))
-\fint_{R(\mathbf{0},0,1)}\ln d_p((y,s),(\mathbf{0},0))
\,dy\,ds\right|\,dx\,dt>0,
\end{align*}
and hence \eqref{20260714.1631} holds. From \eqref{20260714.1631} and
Lemma \ref{Uchiyama}, it follows that $b_0\notin\mathrm{PCMO}$. To show that
$[b_0,I_\alpha^{\gamma+}]_k$ is compact from  $L^r(\omega_0^r)$ to
$L^q(\omega_0^q)$, it suffices to prove that $T$ is compact from $L^r$
to $L^q$, where $T$ is defined by setting, for any suitable function $f$ on
$\mathbb{R}^{n+1}$ and for any $(x,t)\in\mathbb{R}^{n+1}$,
\begin{align*}
Tf(x,t):=\int_{\mathbb{R}^{n+1}}K((x,t),(y,s))f(y,s)\,dy\,ds
\end{align*}
and, for any $(y,s)\in\mathbb{R}^{n+1}$,
\begin{align*}
K((x,t),(y,s)):=\frac{e^{t-s}[b_0(x,t)-b_0(y,s)]^k}
{[d_p((x,t),(y,s))]^{(n+p)(1-\alpha)}}
\boldsymbol{1}_{\Omega_{(x,t)}^{\gamma+}}(y,s).
\end{align*}
We approximate $T$ by a family $\{T_{\varepsilon,M}\}
_{\varepsilon\in(0,1),M\in(\frac{1}{\varepsilon},\infty)}$
of truncation operators  and prove that each
$T_{\varepsilon,M}$ is compact from $L^r$ to $L^q$. To do this, for any
$\varepsilon\in(0,1)$, $M\in(\frac{1}{\varepsilon},\infty)$, and
$(x,t),(y,s)\in\mathbb{R}^{n+1}$, define
\begin{align*}
K_{\varepsilon,M}((x,t),(y,s))&:=K((x,t),(y,s))
\boldsymbol{1}_{\{\varepsilon\leq d_p((x,t),(y,s))\leq\frac{1}{\varepsilon}\}}
((x,t),(y,s))\\
&\quad\times\boldsymbol{1}_{\{d_p((x,t),(\mathbf{0},0))\leq M\}}(x,t)
\boldsymbol{1}_{\{d_p((y,s),(\mathbf{0},0))\leq M\}}(y,s)
\end{align*}
and, for any $f\in L_\mathrm{loc}^1$, let
\begin{align*}
T_{\varepsilon,M}f(x,t):=
\int_{\mathbb{R}^{n+1}}K_{\varepsilon,M}((x,t),(y,s))f(y,s)\,dy\,ds.
\end{align*}
Then $K_{\varepsilon,M}\in L_{(x,t)}^qL_{(y,s)}^{r'}$, and an argument similar to
that used in Step 1 of the proof of Lemma \ref{lem:20260715.1703} yields that
$T_{\varepsilon,M}: L^r\to L^q$ is compact. Furthermore, by the mean value
theorem, we find that, for any $(x,t),(y,s)\in\mathbb{R}^{n+1}$,
\begin{align}\label{20261006.2234}
|b_0(x,t)-b_0(y,s)|\leq\frac{d_p((x,t),(y,s))}
{e+\min\{d_p((x,t),(\mathbf{0},0)),\,d_p((y,s),(\mathbf{0},0))\}}
\leq\frac1ed_p((x,t),(y,s)),
\end{align}
which further implies that
\begin{align}\label{20261006.2223}
|K((x,t),(y,s))|&\,\leq\frac{1}{e^k}
\left[d_p((x-y,t-s),(\mathbf{0},0))\right]^{k-(n+p)(1-\alpha)}
e^{-\gamma[d_p((x-y,t-s),(\mathbf{0},0))]^p}\notag\\
&=:\,H(x-y,t-s).
\end{align}
Let $\sigma:=\frac{1}{1-\alpha}$. Then $1+\frac1q=\frac1r+\frac{1}{\sigma}$ and
\begin{align*}
\|H\|_{L^\sigma}^\sigma&=e^{-k\sigma}\int_{\mathbb{R}^{n+1}}
\left[d_p((x,t),(\mathbf{0},0))\right]^{k\sigma-(n+p)}
e^{-\gamma\sigma[d_p((x,t),(\mathbf{0},0))]^p}\,dx\,dt\\
&\sim\int_0^\infty s^{k\sigma-1}e^{-\gamma\sigma s^p}\,ds<\infty.
\end{align*}
Moreover, from \eqref{20261006.2223} and Young's inequality, we infer that,
for any $f\in L^r$,
\begin{align}\label{20261006.2238}
\left\|Tf-T_{\varepsilon,M}f\right\|_{L^q}
&\leq\left\|\left[H\boldsymbol{1}_{\{d_p(\cdot,(\mathbf{0},0))<\varepsilon\}}
\right]*|f|\,\right\|_{L^q}
+\left\|\left[H\boldsymbol{1}_{\{d_p(\cdot,(\mathbf{0},0))
>\frac{1}{\varepsilon}\}}\right]*|f|\,\right\|_{L^q}\notag\\
&\quad+\left\|\int_{\genfrac{}{}{0pt}{}
{\varepsilon\leq d_p(\cdot,(y,s))\leq\frac{1}{\varepsilon}}
{d_p(\cdot,(\mathbf{0},0))>M}}|K((x,t),(y,s))|\,|f(y,s)|
\,dy\,ds\right\|_{L^q}\notag\\
&\quad+\left\|\int_{\genfrac{}{}{0pt}{}
{\varepsilon\leq d_p(\cdot,(y,s))\leq\frac{1}{\varepsilon}}
{d_p((y,s),(\mathbf{0},0))>M}}|K((x,t),(y,s))|\,|f(y,s)|\,dy\,ds\right\|_{L^q}.
\end{align}
Note that, for any $(x,t),(y,s)\in\mathbb{R}^{n+1}$ with
$\varepsilon\leq d_p((x,t),(y,s))\leq\frac{1}{\varepsilon}$ and either
$d_p((x,t),(\mathbf{0},0))>M$ or $d_p((y,s),(\mathbf{0},0))>M$, it holds that
$\min\{d_p((x,t),(\mathbf{0},0)),\,d_p((y,s),(\mathbf{0},0))\}
\geq M-\frac{1}{\varepsilon}$. This, together with \eqref{20261006.2234},
further implies that, for any $(x,t),(y,s)\in\mathbb{R}^{n+1}$,
\begin{align}\label{20261006.2249}
&|K((x,t),(y,s))|\boldsymbol{1}_{\{\varepsilon\leq
d_p((x,t),(y,s))\leq\frac{1}{\varepsilon}\}}((x,t),(y,s))
\left[\boldsymbol{1}_{\{d_p((x,t),(\mathbf{0},0))>M\}}(x,t)
+\boldsymbol{1}_{\{d_p((y,s),(\mathbf{0},0))>M\}}(y,s)\right]\notag\\
&\quad\leq2\left(\frac{\frac{1}{\varepsilon}}
{e+M-\frac{1}{\varepsilon}}\right)^k
\frac{1}{[d_p((x-y,t-s),(\mathbf{0},0))]^{(n+p)(1-\alpha)}}
\boldsymbol{1}_{\{\varepsilon\leq
d_p((x-y,t-s),(\mathbf{0},0))\leq\frac{1}{\varepsilon}\}}(x-y,t-s)\notag\\
&\quad=:2\left(\frac{\frac{1}{\varepsilon}}{e+M-\frac{1}{\varepsilon}}\right)^k
G_\varepsilon(x-y,t-s).
\end{align}
In addition, $\|G_\varepsilon\|_{L^\sigma}^\sigma
\sim\ln\frac{1}{\varepsilon^2}<\infty$. This, together with
\eqref{20261006.2249} and Young's inequality, further implies that
the third and the fourth terms of \eqref{20261006.2238}
can both be controlled by
\begin{align}\label{20261006.2255}
\left(\frac{\frac{1}{\varepsilon}}{e+M-\frac{1}{\varepsilon}}\right)^k
\|G_\varepsilon\|_{L^\sigma}\|f\|_{L^r}
\sim\left(\frac{\frac{1}{\varepsilon}}{e+M-\frac{1}{\varepsilon}}\right)^k
\left(\ln\frac{1}{\varepsilon^2}\right)^\frac{1}{\sigma}\|f\|_{L^r}\to0
\end{align}
as $M\to\infty$.
On the other hand, from Young's inequality, we deduce that the first and the
second terms of \eqref{20261006.2238} can both be controlled by
\begin{align*}
\left[\left\|H\boldsymbol{1}_{\{d_p(\cdot,(\mathbf{0},0))<\varepsilon\}}
\right\|_{L^\sigma}+\left\|H\boldsymbol{1}_{\{d_p(\cdot,(\mathbf{0},0))
>\frac{1}{\varepsilon}\}}\right\|_{L^\sigma}\right]\|f\|_{L^r}.
\end{align*}
Finally, combining this and \eqref{20261006.2255}, letting $M\to\infty$, and
then letting $\varepsilon\to0^+$, we conclude that
$\|T-T_{\varepsilon,M}\|_{L^r\to L^q}\to0$ as $M\to\infty$ and
$\varepsilon\to0^+$, and hence $T$ is compact from $L^r$ to $L^q$ and
$[b_0,I_\alpha^{\gamma+}]_k$ is compact from  $L^r(\omega_0^r)$ to
$L^q(\omega_0^q)$. This completes the proof of the above claim.

The above example indicates that, in the parabolic weighted setting,
the compactness with respect to a single $A_{r,q}^+(\gamma)$ weight
does not characterize PCMO. This is in contrast to the classical Euclidean
theory for commutators of fractional integral operators, where the weighted
compactness of commutators characterizes the CMO regularity of the symbol under
the classical Muckenhoupt-weight assumptions (see, for instance,
\cite[Theorem 1.4]{gwy(cjm-2021)} and \cite[Theorem 1.3]{wy(pams-2018)}).

Furthermore, the above example also shows that the assumption ``$\omega_1\in
A_{r_1,q_1}^-(\gamma)$'' is \emph{natural} in the sense that the conclusion of
Theorem \ref{parabolic extrapolation of compactness} may be false if
$A_{r_1,q_1}^-(\gamma)$ therein is replaced by $A^+_{r_1,q_1}(\gamma)$.
Indeed, if this assumption in Theorem
\ref{parabolic extrapolation of compactness} can be replaced by
``$\omega_1\in A_{r_1,q_1}^+(\gamma)$'', then it follows that
$[b_0,I_\alpha^{\gamma+}] $ is compact from $L^r(\omega^r)$ to
$L^q(\omega^q)$ for any $\omega\in A^+_{r,q}(\gamma)$, in particular from
$L^r$ to $L^q$. This, together with Theorem \ref{thm:20260714.2255}, further
implies that $b_0\in\mathrm{PCMO}$. This contradiction shows that the above
assumption is natural as stated.
\end{remark}

\section{Proof of Theorem \ref{parabolic interpolation of compactness}}
\label{section4}

In this section, we first establish the parabolic weighted
Fr\'echet--Kolmogorov theorem. Applying this, we then prove Theorem
\ref{parabolic interpolation of compactness}. We begin by recalling the
concept of centered parabolic maximal operators
with time lag introduced by Kinnunen and Saari \cite[Definition
4.1]{ks(apde-2016)}. Let $\gamma\in[0,1)$. The \emph{centered parabolic maximal
operator $\mathcal{M}^{\gamma+}$ with time lag} is defined by setting, for
any $f\in L_{\mathrm{loc}}^1$ and $(x,t)\in\mathbb{R}^{n+1}$,
\begin{align*}
\mathcal{M}^{\gamma+}f(x,t):=\sup_{L\in(0,\infty)}
\fint_{R(x,t,L)^+(\gamma)}|f|.
\end{align*}

The following lemma, which is exactly \cite[Theorem 5.4]{ks(apde-2016)},
indicates that the centered parabolic maximal operator with time lag is also
bounded on parabolic weighted Lebesgue spaces.

\begin{lemma}\label{weighted inequality centered}
Let $\gamma\in(0,1)$, $q\in(1,\infty)$, and $\omega$ be a weight.
Then $\mathcal{M}^{\gamma+}$ is bounded on $L^q(\omega)$ if and only if
$\omega\in A_q^+(\gamma)$.
\end{lemma}

Now, we establish the parabolic weighted Fr\'echet--Kolmogorov theorem, which
characterizes relatively compact subsets of parabolic weighted Lebesgue spaces.

\begin{proposition}\label{Frechet--Kolmogorov}
Let $\gamma\in(0,1)$, $q\in(1,\infty)$, and $\omega\in
A_q^+(\gamma)$. Then a subset $\mathcal{K}\subset L^q(\omega)$ is
relatively compact if and only if all the following statements hold.
\begin{enumerate}
\item[\rm(i)] $\mathcal{K}$ is bounded, i.e.,
$\sup_{f\in\mathcal{K}}\|f\|_{L^q(\omega)}<\infty$.

\item[\rm(ii)] $\mathcal{K}$ is equitight, i.e.,
$\lim_{L\to\infty}\sup_{f\in\mathcal{K}}
\|f\mathbf{1}_{\mathbb{R}^{n+1}
\setminus R(\mathbf{0},0,L)}\|_{L^q(\omega)}=0.$

\item[\rm(iii)] $\mathcal{K}$ is uniformly equicontinuous, i.e.,
$\lim_{l\to0^+}\sup_{f\in\mathcal{K}}
\|f(\cdot)-f_{R(\cdot,l)^+(\gamma)}\|_{L^q(\omega)}=0$.
\end{enumerate}
\end{proposition}

\begin{proof}
We first show the necessity. Assume that $\mathcal{K}\subset L^q(\omega)$ is
relatively compact. From \cite[p.\,13, Theorem]{Yosida}, we deduce that
$\mathcal{K}$ is totally bounded, i.e., for any $\varepsilon\in(0,1)$, there
exist $N\in\mathbb{N}$ and $\{f_j\}_{j\in\mathbb{N}\cap[1,N]}\subset\mathcal{K}$ such
that, for any $f\in\mathcal{K}$, there exists $k\in\mathbb{N}\cap[1,N]$ satisfying
\begin{align}\label{20240530.2154}
\|f_k-f\|_{L^q(\omega)}<\varepsilon.
\end{align}
In addition, since $\mathcal{K}$ is totally bounded, it follows that
$\mathcal{K}$ is bounded, and hence (i) holds.
Moreover, for any $j\in\mathbb{N}\cap[1,N]$,
there exists $L_j\in(0,\infty)$ such that
\begin{align}\label{20240530.2204}
\left\|f_j\mathbf{1}_{\mathbb{R}^{n+1}\setminus
R(\mathbf{0},0,L_j)}\right\|_{L^q(\omega)}<\varepsilon.
\end{align}
Let $L_0:=\max\{L_j:j\in\mathbb{N}\cap[1,N]\}$. Using \eqref{20240530.2154} and
\eqref{20240530.2204}, we find that, for any $f\in\mathcal{K}$ and
$L\in[L_0,\infty)$,
\begin{align*}
\left\|f\mathbf{1}_{\mathbb{R}^{n+1}\setminus R(\mathbf{0},0,L)}\right\|
_{L^q(\omega)}\leq\|f-f_k\|_{L^q(\omega)}+
\left\|f_k\mathbf{1}_{\mathbb{R}^{n+1}\setminus R(\mathbf{0},0,L_k)}\right\|
_{L^q(\omega)}<2\varepsilon.
\end{align*}
Taking the supremum over all $f\in\mathcal{K}$
and letting $\varepsilon\to0^+$, we obtain (ii).
Furthermore, for any $l\in(0,\infty)$ and $f\in\mathcal{K}$,
\begin{align}\label{20240530.2217}
\left\|f(\cdot)-f_{R(\cdot,l)^+(\gamma)}\right\|_{L^q(\omega)}
&\leq\|f-f_k\|_{L^q(\omega)}+\left\|f_k(\cdot)-
(f_k)_{R(\cdot,l)^+(\gamma)}\right\|_{L^q(\omega)}\notag\\
&\quad+\left\|(f_k)_{R(\cdot,l)^+(\gamma)}-
f_{R(\cdot,l)^+(\gamma)}\right\|_{L^q(\omega)}.
\end{align}
From Lemma \ref{weighted inequality centered}, it follows that, for
any $(x,t)\in\mathbb{R}^{n+1}$,
\begin{align}\label{20240530.2232}
\left|f_k(x,t)-(f_k)_{R(x,t,l)^+(\gamma)}\right|
\leq\left|f_k(x,t)\right|+\mathcal{M}^{\gamma+}f_k(x,t)\in L^q(\omega).
\end{align}
In addition, by the Lebesgue differentiation theorem (see
\cite[Lemma 2.3]{kmy(ma-2023)}), we find that, for almost every
$(x,t)\in\mathbb{R}^{n+1}$,
\begin{align*}
(f_k)_{R(x,t,l)^+(\gamma)}\to f_k(x,t)\mbox{\ \ as\ \ }l\to0^+.
\end{align*}
This, together with \eqref{20240530.2232} and the Lebesgue dominated
convergence theorem, further implies that there exists $l_0\in(0,\infty)$ such
that, for any $l\in(0,l_0)$ and $j\in\mathbb{N}\cap[1,N]$,
\begin{align}\label{20240530.2238}
\left\|f_j(\cdot)-(f_j)_{R(\cdot,l)^+(\gamma)}\right\|
_{L^q(\omega)}<\varepsilon.
\end{align}
Note that, for any $(x,t)\in\mathbb{R}^{n+1}$,
\begin{align*}
\left|(f_k)_{R(x,t,l)^+(\gamma)}-f_{R(x,t,l)^{+}(\gamma)}\right|
\leq\fint_{R(x,t,l)^+(\gamma)}|f_k-f|\leq
\mathcal{M}^{\gamma+}(f_k-f)(x,t),
\end{align*}
which, together with Lemma \ref{weighted inequality centered},
further implies that
\begin{align*}
\left\|(f_k)_{R(\cdot,l)^+(\gamma)}-
f_{R(\cdot,l)^+(\gamma)}\right\|
_{L^q(\omega)}\leq\left\|\mathcal{M}^{\gamma+}(f_k-f)\right\|
_{L^q(\omega)}\lesssim\|f_k-f\|_{L^q(\omega)}\leq\varepsilon.
\end{align*}
Combining this, \eqref{20240530.2217}, \eqref{20240530.2238},
and \eqref{20240530.2154} and taking the supremum over all
$f\in\mathcal{K}$, we conclude that, for any $l\in(0,l_0)$,
\begin{align*}
\sup_{f\in\mathcal{K}}\left\|f(\cdot)-f_{R(\cdot,l)^+(\gamma)}\right\|
_{L^q(\omega)}\lesssim\varepsilon,
\end{align*}
and hence (iii) holds. This completes the proof of the necessity.

Next, we prove the sufficiency. Assume that $\mathcal{K}\subset L^q(\omega)$
satisfies (i)--(iii). Then, for any $\varepsilon\in(0,\infty)$, there exist
$0<l_0<L<\infty$ such that, for any $l\in(0,l_0)$ and $f\in\mathcal{K}$,
\begin{align}\label{20240531.1024}
\left\|f\mathbf{1}_{\mathbb{R}^{n+1}
\setminus R(\mathbf{0},0,L)}\right\|_{L^q(\omega)}<\varepsilon
\mbox{\ \ and\ \ }
\left\|f(\cdot)-f_{R(\cdot,l)^+(\gamma)}\right\|
_{L^q(\omega)}<\varepsilon.
\end{align}
Fix $l\in(0,l_0)$. For any $(x_1,t_1),(x_2,t_2)\in
\mathbb{R}^{n+1}$, we abbreviate $R(x_1,t_1,l)^+(\gamma)$ and
$R(x_2,t_2,l)^+(\gamma)$, respectively, as $R_1^+$ and $R_2^+$.
From H\"older's inequality, we infer that
\begin{align}\label{20240531.1030}
\left|f_{R_1^+}\right|\leq\left(\fint_{R_1^+}
|f|^q\omega\right)^\frac{1}{q}
\left(\fint_{R_1^+}\omega^{1-q'}\right)^\frac{1}{q'}
\lesssim\frac{\|f\|_{L^q(\omega)}}{l^{n+p}}
\left[\left(\omega^{1-q'}\right)\left(R_1^+\right)\right]
^\frac{1}{q'}
\end{align}
and
\begin{align}\label{20240531.1033}
\left|f_{R_1^+}-f_{R_2^+}\right|
&\leq\frac{1}{|R_1^+|}\int_{\mathbb{R}^{n+1}}
\left|\mathbf{1}_{R_1^+}-
\mathbf{1}_{R_2^+}\right|\,|f|\notag\\
&\lesssim\frac{\|f\|_{L^q(\omega)}}{l^{n+p}}
\left(\int_{\mathbb{R}^{n+1}}
\left|\mathbf{1}_{R_1^+}-
\mathbf{1}_{R_2^+}\right|^{q'}\omega^{1-q'}
\right)^\frac{1}{q'}.
\end{align}
Since $\omega\in A_q^+(\gamma)$, it follows that $\omega^{1-q'}\in
L_{\mathrm{loc}}^1$. Thus, for any $(x,t)\in\overline{R(\mathbf{0},0,L)}$,
\begin{align*}
\left(\omega^{1-q'}\right)\left(R(x,t,l)^+(\gamma)\right)
\leq\left(\omega^{1-q'}\right)\left(\overline{R(\mathbf{0},0
,L+l)}\right)<\infty,
\end{align*}
which, together with \eqref{20240531.1033}, (i), and the Lebesgue
dominated convergence theorem, further implies that
$\{f_{R(\cdot,l)^+(\gamma)}\}_{f\in\mathcal{K}}$ is equicontinuous on
$\overline{R(\mathbf{0},0,L)}$. Combining this,
\eqref{20240531.1030}, and the Arzel\'a--Ascoli theorem (see, for
instance, \cite[p.\,85, Theorem]{Yosida}), we conclude that
$\{f_{R(\cdot,l)^+(\gamma)}\}_{f\in\mathcal{K}}$ is a relatively compact subset
of $C(\overline{R(\mathbf{0},0,L)})$, and hence a totally bounded
subset of $C(\overline{R(\mathbf{0},0,L)})$. Therefore, for any
$\varepsilon\in(0,\infty)$, there exist $M\in\mathbb{N}$ and
$\{f_{j}\}_{j\in\mathbb{N}\cap[1,M]}\subset\mathcal{K}$ such that, for any
$f\in\mathcal{K}$, there exists $k\in\mathbb{N}\cap[1,M]$ satisfying
\begin{align}\label{20240531.1103}
&\left\|f_{R(\cdot,l)^+(\gamma)}-
(f_{k})_{R(\cdot,l)^+(\gamma)}\right\|
_{C(\overline{R(\mathbf{0},0,L)})}\notag\\
&\quad:=\sup_{(x,t)\in\overline{R(\mathbf{0},0,L)}}
\left|f_{R(x,t,l)^+(\gamma)}-
(f_{k})_{R(x,t,l)^+(\gamma)}\right|<\varepsilon
\left[\omega\left(R(\mathbf{0},0,L)\right)+1\right]
^{-\frac{1}{q}}.
\end{align}
To proceed, we split
\begin{align}\label{20240531.1111}
\|f-f_k\|_{L^q(\omega)}\leq\left\|(f-f_k)
\mathbf{1}_{R(\mathbf{0},0,L)}\right\|
_{L^q(\omega)}+\left\|(f-f_k)
\mathbf{1}_{\mathbb{R}^{n+1}\setminus
R(\mathbf{0},0, L)}\right\|_{L^q(\omega)}
=:\mathrm{I}+\mathrm{II}.
\end{align}
We estimate $\mathrm{I}$ and $\mathrm{II}$ separately. On the
one hand, from \eqref{20240531.1024}, we infer that
\begin{align}\label{20240531.1115}
\mathrm{II}\leq\left\|f\mathbf{1}
_{\mathbb{R}^{n+1}\setminus R(\mathbf{0},0,L)}
\right\|_{L^{q}(\omega)}+\left\|f_{k}
\mathbf{1}_{\mathbb{R}^{n+1}\setminus
R(\mathbf{0},0,L)}\right\|_{L^{q}(\omega)}
<2\varepsilon.
\end{align}
On the other hand, by \eqref{20240531.1024} and \eqref{20240531.1103}, we obtain
\begin{align*}
\mathrm{I}&\leq\left\|f(\cdot)-f_{R(\cdot,l)^+(\gamma)}\right\|
_{L^q(\omega)}+\left\|\left[f_{R(\cdot,l)^+(\gamma)}
-(f_k)_{R(\cdot,l)^+(\gamma)}\right]
\mathbf{1}_{R(\mathbf{0},0,L)}\right\|_{L^q(\omega)}\\
&\quad+\left\|(f_k)_{R(\cdot,l)^+(\gamma)}
-f_k(\cdot)\right\|_{L^q(\omega)}\\
&<\varepsilon+\varepsilon\left[\omega\left(
R(\mathbf{0},0,L)\right)+1\right]^{-\frac{1}{q}
}\left\|\mathbf{1}_{R(\mathbf{0},0,L)}\right\|
_{L^q(\omega)}+\varepsilon<3\varepsilon.
\end{align*}
This, together with \eqref{20240531.1111} and
\eqref{20240531.1115}, further implies that, for any
$\varepsilon\in(0,\infty)$, there exist $M\in\mathbb{N}$ and
$\{f_j\}_{j\in\mathbb{N}\cap[1,M]}\subset\mathcal{K}$ such that, for any $f\in
\mathcal{K}$, there exists $k\in\mathbb{N}\cap[1,M]$ satisfying
\begin{align*}
\|f-f_k\|_{L^q(\omega)}<5\varepsilon.
\end{align*}
Therefore, $\mathcal{K}$ is totally bounded, and hence is relatively
compact in $L^q(\omega)$. This completes the proof of
the sufficiency and hence Proposition \ref{Frechet--Kolmogorov}.
\end{proof}

Proposition \ref{Frechet--Kolmogorov} further implies the following
characterization of relatively compact subsets of parabolic weighted
Lebesgue spaces $L^q(\omega)$ with $q\in(0,\infty)$ and
$\omega\in A_\infty^+(\gamma)$, which plays a crucial role in the proof of
Theorem \ref{parabolic interpolation of compactness}.

\begin{proposition}\label{Frechet--Kolmogorov 2}
Let $\gamma\in(0,1)$, $q\in(0,\infty)$, and $\omega\in A_\infty^+(\gamma)$.
Let $q_0\in(1,\infty)$ be such that $\omega\in A_{q_0}^+(\gamma)$.
Then a subset $\mathcal{K}\subset L^q(\omega)$ is relatively compact
if and only if all the following assertions hold.
\begin{enumerate}
\item[\rm(i)] $\mathcal{K}$ is bounded, i.e.,
$\sup_{f\in\mathcal{K}}\|f\|_{L^q(\omega)}<\infty$.

\item[\rm(ii)] $\mathcal{K}$ is equitight, i.e.,
$\lim_{L\to\infty}\sup_{f\in\mathcal{K}}
\|f\boldsymbol{1}_{\mathbb{R}^{n+1}\setminus R(\mathbf{0},0,L)}
\|_{L^q(\omega)}=0$.

\item[\rm(iii)]
\begin{align*}
\lim_{l\to0^+}\sup_{f\in\mathcal{K}}
\int_{\mathbb{R}^{n+1}}\left[\fint_{R(x,t,l)^+(\gamma)}
|f(x,t)-f(y,s)|^\frac{q}{q_0}\,dy\,ds\right]^{q_0}\omega(x,t)\,dx\,dt=0.
\end{align*}
\end{enumerate}
\end{proposition}

\begin{proof}
We first show the necessity. Assume that $\mathcal{K}\subset L^q(\omega)$ is
relatively compact. Observe that, for any $f,g\in L^q(\omega)$,
\begin{align*}
\|f+g\|_{L^q(\omega)}\leq\max\left\{1,\,2^{\frac1q-1}\right\}
\left[\|f\|_{L^q(\omega)}+\|g\|_{L^q(\omega)}\right].
\end{align*}
Thus, by an argument similar to that used in the proof of Proposition
\ref{Frechet--Kolmogorov}, we conclude that both (i) and (ii)
hold. Then we justify (iii). Since $\mathcal{K}$ is totally bounded, it follows
that, for any $\varepsilon\in(0,\infty)$, there exist
$N\in\mathbb{N}$ and $\{f_j\}_{j\in\mathbb{N}\cap[1,N]}\subset\mathcal{K}$ such
that, for any $f\in\mathcal{K}$, there exists $k\in\mathbb{N}\cap[1,N]$
satisfying \eqref{20240530.2154}. Note that, for any $l\in(0,\infty)$,
\begin{align}\label{20140601.1505}
&\int_{\mathbb{R}^{n+1}}\left[\fint_{R(x,t,l)^+(\gamma)}
|f(x,t)-f(y,s)|^\frac{q}{q_0}\,dy\,ds\right]^{q_0}\omega(x,t)\,dx\,dt\notag\\
&\quad\lesssim\int_{\mathbb{R}^{n+1}}
\left[\fint_{R(x,t,l)^+(\gamma)}
|f(x,t)-f_k(x,t)|^\frac{q}{q_0}\,dy\,ds\right]^{q_0}\omega(x,t)\,dx\,dt\notag\\
&\quad\quad+\int_{\mathbb{R}^{n+1}}
\left[\fint_{R(x,t,l)^+(\gamma)}|f_k(x,t)-f_k(y,s)|
^\frac{q}{q_0}\,dy\,ds\right]^{q_0}\omega(x,t)\,dx\,dt\notag\\
&\quad\quad+\int_{\mathbb{R}^{n+1}}
\left[\fint_{R(x,t,l)^+(\gamma)}|f_k(y,s)-f(y,s)|
^\frac{q}{q_0}\,dy\,ds\right]^{q_0}\omega(x,t)\,dx\,dt\notag\\
&\quad=:\mathrm{I}+\mathrm{II}+\mathrm{III}.
\end{align}
To estimate $\mathrm{I}$, using
\eqref{20240530.2154}, we obtain
\begin{align}\label{20240601.1521}
\mathrm{I}=\int_{\mathbb{R}^{n+1}}|f-f_k|^q\omega<\varepsilon^q.
\end{align}
To estimate $\mathrm{III}$, by Lemma \ref{weighted inequality centered}
and \eqref{20240530.2154}, we conclude that
\begin{align}\label{20240601.1526}
\mathrm{III}\leq\int_{\mathbb{R}^{n+1}}
\left[\mathcal{M}^{\gamma+}\left(|f-f_k|^\frac{q}{q_0}\right)^{q_0}\right]
\omega\lesssim\int_{\mathbb{R}^{n+1}}|f-f_k|^q\omega<\varepsilon^q.
\end{align}
We turn to estimate $\mathrm{II}$. Denote by the \emph{notation $C_\mathrm{c}$} the
set of all continuous functions on $\mathbb{R}^{n+1}$ with compact support.
Since $C_\mathrm{c}$ is dense in $L^q(\omega)$, it follows that, for any
$k\in\mathbb{N}\cap[1,N]$, there exists $g_k\in C_\mathrm{c}$ such that
\begin{align}\label{20240601.1536}
\|f_k-g_k\|_{L^q(\omega)}<\varepsilon.
\end{align}
Moreover, there exist $0<l_0<1<L_0<\infty$ such that, for any
$k\in\mathbb{N}\cap[1,N]$, $\mathrm{supp}(g_k)\subset R(\mathbf{0},0,L_0)$
and, for any $(x,t),(y,s)\in\mathbb{R}^{n+1}$ with $d_p((x,t),(y,s))\in[0,l_0]$,
\begin{align}\label{20240601.1540}
|g_k(x,t)-g_k(y,s)|<\varepsilon
\left[1+\omega(R(\mathbf{0},0,L_0+1))\right]^{-\frac1q}.
\end{align}
Applying \eqref{20240601.1536}, \eqref{20240601.1540}, and Lemma
\ref{weighted inequality centered}, we conclude that, for any $l\in(0,l_0)$,
\begin{align*}
\mathrm{II}&\lesssim\int_{\mathbb{R}^{n+1}}|f_k(x,t)-g_k(x,t)|^q
\omega(x,t)\,dx\,dt\\
&\quad+\int_{\mathbb{R}^{n+1}}\left[\fint_{R(x,t,l)^+(\gamma)}
|g_k(x,t)-g_k(y,s)|^\frac{q}{q_0}\,dy\,ds\right]^{q_0}\omega(x,t)\,dx\,dt\\
&\quad\quad+\int_{\mathbb{R}^{n+1}}\left[\fint_{R(x,t,l)^+(\gamma)}
|g_k(y,s)-f_k(y,s)|^\frac{q}{q_0}\,dy\,ds\right]^{q_0}\omega(x,t)\,dx\,dt\\
&\lesssim\int_{\mathbb{R}^{n+1}}|f_k-g_k|^q\omega
+\varepsilon^q\left[1+\omega(R(\mathbf{0},0,L_0+1))\right]^{-1}
\omega(R(\mathbf{0},0,L_0+1))\\
&\quad+\int_{\mathbb{R}^{n+1}}\mathcal{M}^{\gamma+}
\left(|g_k-f_k|^\frac{q}{q_0}\right)^{q_0}\omega\\
&\lesssim\int_{\mathbb{R}^{n+1}}|f_k-g_k|^q\omega+\varepsilon^q
\lesssim\varepsilon^q.
\end{align*}
From this, \eqref{20140601.1505}, \eqref{20240601.1521}, and
\eqref{20240601.1526} and taking the supremum over all $f\in\mathcal{K}$, we
deduce that, for any $l\in(0,l_0)$,
\begin{align*}
\sup_{f\in\mathcal{K}}\int_{\mathbb{R}^{n+1}}
\left[\fint_{R(x,t,l)^+(\gamma)}|f(x,t)-f(y,s)|^\frac{q}{q_0}
\,dy\,ds\right]^{q_0}\omega(x,t)\,dx\,dt\lesssim\varepsilon^q.
\end{align*}
This completes the proof of (iii) and hence the necessity.

Now, we prove the sufficiency. Suppose that $\mathcal{K}\subset L^q(\omega)$
satisfies (i)--(iii). We consider the following two cases for $q$ and $q_0$.

\emph{Case (1)} $1<q_0\leq q<\infty$. In this case, from H\"older's inequality,
we infer that, for any $f\in\mathcal{K}$ and $(x,t)\in\mathbb{R}^{n+1}$,
\begin{align*}
\left|f(x,t)-f_{R(x,t,l)^+(\gamma)}\right|
&\leq\fint_{R(x,t,l)^+(\gamma)}|f(x,t)-f(y,s)|\,dy\,ds\\
&\leq\left[\fint_{R(x,t,l)^+(\gamma)}
|f(x,t)-f(y,s)|^\frac{q}{q_0}\,dy \,ds\right]^\frac{q_0}{q},
\end{align*}
which, together with (iii), further implies that
\begin{align}\label{20240601.1640}
\lim_{l\to0^+}\sup_{f\in\mathcal{K}}
\left\|f(\cdot)-f_{R(\cdot,l)^+(\gamma)}\right\|_{L^q(\omega)}=0.
\end{align}
In addition, since $1<q_0\leq q<\infty$ and $\omega\in
A_{q_0}^+(\gamma)$,  it follows that $\omega\in A_q^+(\gamma)$.
Applying this, (i), (ii), \eqref{20240601.1640}, and Proposition
\ref{Frechet--Kolmogorov}, we find that $\mathcal{K}$ is relatively compact.

\emph{Case (2)} $0<q<q_{0}<\infty$. In this case, we divide the
proof into the following three steps.

\emph{Step 1.} Assume that all the functions in $\mathcal{K}$ are
non-negative. By a simple calculation, we obtain the
following basic inequality: for any $\alpha\in(0,1]$ and $a,b\in
[0,\infty)$ with $a\neq b$,
\begin{align}\label{20250115.1935}
\left|a^\alpha-b^\alpha\right|\leq|a-b|^\alpha\leq
\frac{1}{\alpha}\left(\frac{a+b}{|a-b|}\right)^{1-\alpha}
\left|a^\alpha-b^\alpha\right|
\end{align}
(see \cite[p.\,33]{t(kmsr-1951)}). From the first inequality of
\eqref{20250115.1935} with $\alpha:=\frac{q}{q_0}$ therein and from H\"older's
inequality, we deduce that, for any $f\in\mathcal{K}$ and
$(x,t)\in\mathbb{R}^{n+1}$,
\begin{align*}
\left|\left[f(x,t)\right]^\frac{q}{q_0}-
\left(f^\frac{q}{q_0}\right)_{R(x,t,l)^+(\gamma)}\right|
\leq\fint_{R(x,t,l)^+(\gamma)}\left|f(x,t)-f(y,s)\right|^\frac{q}{q_0}\,dy\,ds,
\end{align*}
which, combined with (iii), further implies that
\begin{align}\label{20240601.2226}
\lim_{l\to0^+}\sup_{f\in\mathcal{K}}\left\|\left[f(\cdot)\right]^\frac{q}{q_0}
-\left(f^\frac{q}{q_0}\right)_{R(\cdot,l)^+(\gamma)}
\right\|_{L^{q_0}(\omega)}=0.
\end{align}
Moreover, the conditions (i) and (ii) can be
rewritten, respectively, as
\begin{align*}
\sup_{f\in\mathcal{K}}\left\|f^\frac{q}{q_0}\right\|_{L^{q_0}(\omega)}
<\infty\mbox{\ \ and\ \ }\lim_{L\to\infty}\sup_{f\in\mathcal{K}}
\left\|f^\frac{q}{q_0}\boldsymbol{1}_{\mathbb{R}^{n+1}
\setminus R(\mathbf{0},0,L)}\right\|_{L^{q_0}(\omega)}=0.
\end{align*}
This, together with \eqref{20240601.2226} and Proposition
\ref{Frechet--Kolmogorov}, further implies that
$\mathcal{K}^\frac{q}{q_0}:=\{f^\frac{q}{q_0}:f\in\mathcal{K}\}$
is relatively compact in $L^{q_0}(\omega)$.

\emph{Step 2.} We still assume that all the functions in $\mathcal{K}$ are
non-negative. Since $\mathcal{K}^\frac{q}{q_0}$ is relatively compact in
$L^{q_0}(\omega)$, it follows that any sequence
$\{f_j\}_{j\in\mathbb{N}}\subset\mathcal{K}$ has a subsequence, still
denoted by $\{f_j\}_{j\in\mathbb{N}}$ for simplicity of presentation, such that
$\{f_j^\frac{q}{q_0}\}_{j\in\mathbb{N}}$ is a Cauchy sequence in
$L^{q_0}(\omega)$. Therefore, for any $\varepsilon\in(0,\infty)$, there
exists $N\in\mathbb{N}$ such that, for all $i,j\in\mathbb{N}\cap[N,\infty)$,
\begin{align}\label{20240601.2243}
\int_{\mathbb{R}^{n+1}}\left|f_i^\frac{q}{q_0}-
f_j^\frac{q}{q_0}\right|^{q_0}\omega<\varepsilon^{q_0}.
\end{align}
For any given $i,j\in\mathbb{N}\cap[N,\infty)$, define
\begin{align}\label{20240601.2303}
E_\varepsilon^{i,j}:=\left\{(x,t)\in\mathbb{R}^{n+1}:
0\leq f_i(x,t)+f_j(x,t)\leq\frac{|f_i(x,t)-f_j(x,t)|}{\varepsilon}\right\}.
\end{align}
Applying \eqref{20240601.2303}, the second inequality of
\eqref{20250115.1935} with $\alpha:=\frac{q}{q_0}$ therein, and
\eqref{20240601.2243}, we find that, for any $i,j\in\mathbb{N}\cap[N,\infty)$,
\begin{align}\label{20240601.2258}
\int_{E_\varepsilon^{i,j}}\left|f_i-f_j\right|^q\omega
\leq\left(\frac{q_0}{q}\right)^{q_0}\varepsilon^{(\frac{q}{q_{0}}-1)q_0}
\int_{E_\varepsilon^{i,j}}\left|f_i^\frac{q}{q_{0}}-
f_j^\frac{q}{q_{0}}\right|^{q_0}\omega
\leq\left(\frac{q_{0}}{q}\right)^{q_0}\varepsilon^q.
\end{align}
In addition, by \eqref{20240601.2303} and (i), we conclude
that, for any $i,j\in\mathbb{N}\cap[N,\infty)$,
\begin{align*}
\int_{\mathbb{R}^{n+1}\setminus E_\varepsilon^{i,j}}
\left|f_i-f_j\right|^q\omega
&\leq \int_{\mathbb{R}^{n+1}\setminus E_\varepsilon^{i,j}}\varepsilon^q
\left|f_i+f_j\right|^q\omega\\
&\lesssim\varepsilon^q
\left(\int_{\mathbb{R}^{n+1}\setminus E_\varepsilon^{i,j}}|f_i|^q\omega
+\int_{\mathbb{R}^{n+1}\setminus E_\varepsilon^{i,j}}
\left|f_j\right|^q\omega\right)\\
&\lesssim\sup_{f\in\mathcal{K}}
\|f\|_{L^q(\omega)}^q\varepsilon^q.
\end{align*}
Combining this and \eqref{20240601.2258}, we find that the selected
subsequence $\{f_j\}_{j\in\mathbb{N}}$ is a Cauchy sequence in $L^q(\omega)$.
Since $L^q(\omega)$ is complete, this subsequence converges in
$L^q(\omega)$. Thus, any sequence in $\mathcal{K}$ has a convergent subsequence,
and hence $\mathcal{K}$ is relatively compact in $L^q(\omega)$.

\emph{Step 3.} We now consider the general case. Define
\begin{align*}
\mathcal{K}_+:=\{f_+:f\in\mathcal{K}\}\mbox{\ \ and\ \ }
\mathcal{K}_-:=\{f_-:f\in\mathcal{K}\}.
\end{align*}
Here, and thereafter, for any function $f$ on $\mathbb{R}^{n+1}$, let
$f_+:=\max\{f,\,0\}$ and $f_-=-\min\{f,\,0\}$.
For any $f\in\mathcal{K}$ and $(x,t),(y,s)\in\mathbb{R}^{n+1}$,
$0\leq f_\pm(x,t)\leq |f(x,t)|$ and
$|f_\pm(x,t)-f_\pm(y,s)|\leq|f(x,t)-f(y,s)|$.
Thus, both $\mathcal{K}_+$ and $\mathcal{K}_-$ satisfy (i)--(iii), and hence are
relatively compact in $L^q(\omega)$ by Step 2. Given a sequence
$\{f_j\}_{j\in\mathbb{N}}\subset\mathcal{K}$, we first select a subsequence
whose positive parts converge and then a further subsequence whose negative
parts converge in $L^q(\omega)$. Since $f_j=(f_j)_+-(f_j)_-$,
the quasi-triangle inequality yields that the resulting subsequence converges
in $L^q(\omega)$. Consequently, any sequence in $\mathcal{K}$ has a convergent
subsequence in $L^q(\omega)$, and hence $\mathcal{K}$ is relatively compact in
$L^q(\omega)$. This completes the proof of the sufficiency and hence of
Proposition \ref{Frechet--Kolmogorov 2}.
\end{proof}

Finally, we are ready to show Theorem
\ref{parabolic interpolation of compactness}.

\begin{proof}[Proof of Theorem \ref{parabolic interpolation of compactness}]
Fix $\theta\in(0,1)$. To show that $T$ is compact from
$L^{r_\theta}(\omega_\theta^{r_\theta})$ to
$L^{q_\theta}(\omega_\theta^{q_\theta})$, it suffices to prove that $
T(\overline{B_\theta})$ is relatively compact in
$L^{q_\theta}(\omega_\theta^{q_\theta})$, where
\begin{align*}
\overline{B_\theta}:=
\left\{f\in L^{r_\theta}(\omega_\theta^{r_\theta}):
\|f\|_{L^{r_\theta}(\omega_\theta^{r_\theta})}\leq1\right\}.
\end{align*}
From the definition of $A_\infty^+(\gamma)$, it follows that there
exists $r\in(1,\infty)$ such that $\omega_0^{q_0},
\omega_1^{q_1}\in A_r^+(\gamma)$. Applying the assumptions that $T$ is bounded
from $L^{r_0}(\omega_0^{r_0})$ to $L^{q_0}(\omega_0^{q_0})$ and
is compact from $L^{r_1}(\omega_1^{r_1})$ to $L^{q_1}(\omega_1^{q_1})$, we
conclude that there exist positive constants $A_0$ and $A_1$ such that, for any
$f\in L^{r_0}(\omega_0^{r_0})$,
\begin{align}\label{20240602.2120}
\|Tf\|_{L^{q_0}(\omega_0^{q_0})}\leq A_0\|f\|_{L^{r_0}(\omega_0^{r_0})}
\end{align}
and, for any $f\in L^{r_1}(\omega_1^{r_1})$,
\begin{align*}
\|Tf\|_{L^{q_1}(\omega_1^{q_1})}\leq A_1\|f\|_{L^{r_1}(\omega_1^{r_1})}.
\end{align*}
From this, \eqref{20240602.2120}, and the interpolation theorem (see
\cite[Lemmas 3.2 and 3.3]{coy(tams-2022)}), we deduce that,
for any $f\in\overline{B_\theta}$,
\begin{align}\label{20240602.2146}
\|Tf\|_{L^{q_{\theta}}(\omega_{\theta}^{q_\theta})}\leq
A_0^{1-\theta}A_1^{\theta}.
\end{align}
Moreover, using the assumption that $T$ is compact from
$L^{r_1}(\omega_1^{r_1})$ to $L^{q_1}(\omega_1^{q_1})$ and applying
Proposition \ref{Frechet--Kolmogorov 2} to $T(\mathcal{K}_1)$ with
$\mathcal{K}_1:=\{f\in L^{r_1}(\omega_1^{r_1}):
\|f\|_{L^{r_1}(\omega_1^{r_1})}\leq1\}$, we obtain
\begin{align}\label{20240602.2129}
\lim_{L\to\infty}\sup_{\|f\|_{L^{r_1}(\omega_1^{r_1})}\leq1}
\left\|Tf\boldsymbol{1}_{\mathbb{R}^{n+1}
\setminus R(\mathbf{0},0,L)}\right\|_{L^{q_1}(\omega_1^{q_1})}=0
\end{align}
and
\begin{align}\label{20240602.2131}
\lim_{l\to0^+}\sup_{\|f\|_{L^{r_1}(\omega_1^{r_1})}\leq1}
\int_{\mathbb{R}^{n+1}}\left[\fint_{R(x,t,l)^+(\gamma)}
\left|Tf(x,t)-Tf(y,s)\right|^\frac{q_1}{r}\,dy\,ds\right]^r
[\omega_1(x,t)]^{q_1}\,dx\,dt=0.
\end{align}
From \eqref{20240602.2129}, we infer that, for any
$\varepsilon\in(0,\infty)$, there exists $L_\varepsilon\in(0,\infty)$
such that, for any $L\in(L_\varepsilon,\infty)$ and $f\in
L^{r_1}(\omega_1^{r_1})$ with $\|f\|_{L^{r_1}(\omega_1^{r_1})}
\in(0,\infty)$,
\begin{align}\label{20240602.2156}
\left\|Tf\mathbf{1}_{\mathbb{R}^{n+1}\setminus
R(\mathbf{0},0,L)}\right\|_{L^{q_1}(\omega_1^{q_1})}\leq\varepsilon
\|f\|_{L^{r_1}(\omega_1^{r_1})}.
\end{align}
On the other hand, \eqref{20240602.2156} obviously holds for any
$f\in L^{r_1}(\omega_1^{r_1})$ with $\|f\|_{L^{r_1}(\omega_1^{r_1})}=0$, and
hence holds for all $f\in L^{r_1}(\omega_1^{r_1})$. This, together with
\eqref{20240602.2120} and \cite[Lemmas 3.2 and 3.3]{coy(tams-2022)}, further
implies that, for any $L\in(L_\varepsilon,\infty)$ and $f\in
L^{r_\theta}(\omega_\theta^{r_\theta})$,
\begin{align*}
\left\|Tf\mathbf{1}_{\mathbb{R}^{n+1}\setminus
R(\mathbf{0},0,L)}\right\|_{L^{q_\theta}(\omega_\theta^{q_\theta})}\leq
A_0^{1-\theta}\varepsilon^{\theta}
\|f\|_{L^{r_\theta}(\omega_\theta^{r_\theta})}.
\end{align*}
Thus,
\begin{align}\label{20240602.2203}
\lim_{L\to\infty}\sup_{f\in \overline{B_\theta}}
\left\|Tf\boldsymbol{1}_{\mathbb{R}^{n+1}\setminus
R(\mathbf{0},0,L)}\right\|_{L^{q_\theta}
(\omega_\theta^{q_\theta})}=0.
\end{align}
By \eqref{20240602.2131}, we find that there exists
$l_\varepsilon\in(0,\infty)$ such that, for any
$l\in(0,l_\varepsilon)$ and $f\in L^{r_1}(\omega_1^{r_1})$ with
$\|f\|_{L^{r_1}(\omega_1^{r_1})}\in(0,\infty)$,
\begin{align}\label{20240602.2213}
\int_{\mathbb{R}^{n+1}}\left[\fint_{R(x,t,l)^+(\gamma)}
|Tf(x,t)-Tf(y,s)|^\frac{q_1}{r}\,dy\,ds\right]^r
[\omega_1(x,t)]^{q_1}\,dx\,dt
\leq\varepsilon^{q_1}\|f\|_{L^{r_1}(\omega_1^{r_1})}^{q_1}.
\end{align}
On the other hand, \eqref{20240602.2213} holds automatically for
any $f\in L^{r_1}(\omega_1^{r_1})$ with
$\|f\|_{L^{r_1}(\omega_1^{r_1})}=0$, and hence holds for all
$f\in L^{r_1}(\omega_1^{r_1})$. Moreover, from Lemma
\ref{weighted inequality centered} and \eqref{20240602.2120}, we deduce that,
for any $l\in(0,\infty)$ and $f\in L^{r_0}(\omega_0^{r_0})$,
\begin{align*}
&\int_{\mathbb{R}^{n+1}}\left[
\fint_{R(x,t,l)^+(\gamma)} |Tf(x,t)-Tf(y,s)|^\frac{q_0}{r}
\,dy\,ds\right]^r[\omega_0(x,t)]^{q_0}\,dx\,dt\\
&\quad\lesssim\int_{\mathbb{R}^{n+1}}|Tf|^{q_0}\omega_0^{q_0}
+\int_{\mathbb{R}^{n+1}}\left[\mathcal{M}^{\gamma+}\left(|Tf|
^\frac{q_0}{r}\right)\right]^r\omega_0^{q_0}\\
&\quad\lesssim\|Tf\|_{L^{q_0}(\omega_0^{q_0})}^{q_0}
\leq A_0^{q_0}\|f\|_{L^{r_0}(\omega_0^{r_0})}^{q_0}.
\end{align*}
This, along with \eqref{20240602.2213} and the
weighted interpolation on mixed-norm Lebesgue spaces
(see \cite[Theorem 3.5]{coy(tams-2022)}),
further implies that, for any
$l\in(0,l_\varepsilon)$ and $f\in
L^{r_\theta}(\omega_\theta^{r_\theta})$,
\begin{align*}
&\int_{\mathbb{R}^{n+1}}\left[\fint_{R(x,t,l)^+(\gamma)}
|Tf(x,t)-Tf(y,s)|^\frac{q_\theta}{r}\,dy \,ds\right]^r
[\omega_\theta(x,t)]^{q_\theta}\,dx\,dt\\
&\quad\lesssim A_0^{(1-\theta)q_\theta}\varepsilon^{\theta q_\theta}
\|f\|_{L^{r_\theta}(\omega_\theta^{r_\theta})}^{q_\theta},
\end{align*}
which means that
\begin{align}\label{1520}
\lim_{l\to0^+}\sup_{f\in\overline{B_\theta}}
\int_{\mathbb{R}^{n+1}}\left[\fint_{R(x,t,l)^+(\gamma)}
|Tf(x,t)-Tf(y,s)|^\frac{q_\theta}{r}\,dy \,ds\right]^r
[\omega_\theta(x,t)]^{q_\theta}\,dx\,dt=0.
\end{align}

To apply Proposition \ref{Frechet--Kolmogorov 2} to the set
$T(\overline{B_\theta})$, we still need to verify
$\omega_\theta^{q_\theta}\in A_\infty^+(\gamma)$.
Indeed, from $\omega_\theta=\omega_0^{1-\theta}\omega_1^\theta$,
$\frac{1}{q_\theta}=\frac{(1-\theta)}{q_0}+\frac{\theta}{q_1}$,
and H\"older's inequality twice, we deduce that
\begin{align*}
\left[\omega_\theta^{q_\theta}\right]_{A_r^+(\gamma)}
&=\sup_{R\in\mathcal{R}}\fint_{R^-(\gamma)}
\left(\omega_0^{q_0}\right)^{\frac{(1-\theta)q_\theta}{q_0}}
\left(\omega_1^{q_1}\right)^{\frac{\theta q_\theta}{q_1}}\\
&\quad\times\left[\fint_{R^+(\gamma)}
\left(\omega_0^{q_0}\right)^{\frac{(1-r')(1-\theta)q_\theta}{q_0}}
\left(\omega_1^{q_1}\right)^{\frac{(1-r')\theta q_\theta}{q_1}}\right]^{r-1}\\
&\leq\left[\omega_0^{q_0}\right]_{A_r^+(\gamma)}^\frac{(1-\theta)q_\theta}{q_0}
\left[\omega_1^{q_1}\right]_{A_r^+(\gamma)}^\frac{\theta q_\theta}{q_1}<\infty.
\end{align*}
Therefore, $\omega_\theta^{q_\theta}\in A_r^+(\gamma)
\subset A_\infty^+(\gamma)$.
Combining this, \eqref{20240602.2203},
\eqref{20240602.2146}, \eqref{1520},
and Proposition \ref{Frechet--Kolmogorov 2},
we conclude that $T(\overline{B_\theta})$ is relatively compact
in $L^{q_\theta}(\omega_\theta^{q_\theta})$. This completes
the proof of Theorem
\ref{parabolic interpolation of compactness}.
\end{proof}

\medskip

\noindent\textbf{Acknowledgements}\quad
The main mathematical ideas and approaches in this article
originated from the authors. Generative AI tools were used as
supplementary aids during the exploratory and manuscript-preparation
stages, particularly in refining the presentation and examining
technical details. All mathematical arguments and proofs in the
final manuscript were independently verified and finalized by
the authors, who take full responsibility for their correctness.

\bigskip

\noindent Mingming Cao

\medskip

\noindent Institute for Advanced Study in Mathematics,
Harbin Institute of Technology, Harbin 150001, The People's Republic of China

\smallskip

\noindent{\it E-mail:} \texttt{mingming.cao@hit.edu.cn}

\bigskip

\noindent Weiyi Kong, Dachun Yang and Wen Yuan

\medskip

\noindent Laboratory of Mathematics and Complex Systems (Ministry of Education of China),
School of Mathematical Sciences, Institute for Advanced Study,
Beijing Normal University, Beijing 100875,
The People's Republic of China

\smallskip

\noindent{\it E-mails:} \texttt{weiyikong@mail.bnu.edu.cn} (W. Kong)

\noindent\phantom{{\it E-mails:} }\texttt{dcyang@bnu.edu.cn} (D. Yang)

\noindent\phantom{{\it E-mails:} }\texttt{wenyuan@bnu.edu.cn} (W. Yuan)

\bigskip

\noindent Chenfeng Zhu

\medskip

\noindent School of Mathematical Sciences, Zhejiang
University of Technology, Hangzhou 310023, The People's Republic of China

\smallskip

\noindent{\it E-mail:} \texttt{chenfengzhu@zjut.edu.cn}

\end{document}